\documentclass[a4paper]{amsart}

\title[Vertex algebras associated with hypertoric varieties]{%
Vertex algebras associated with hypertoric varieties}
\author{Toshiro Kuwabara}
\thanks{The author is partially supported by Grant-in-Aid for Young Scientist (B) 17K14151,
Japan Society for the Promotion of Science}

\usepackage{latexsym}
\usepackage{amsmath}
\usepackage{amsfonts}
\usepackage{amssymb}
\usepackage{amsxtra}
\usepackage{mathrsfs}

\newtheorem{definition}{Definition}[section]
\newtheorem{proposition}[definition]{Proposition}

\newtheorem{lemma}[definition]{Lemma}
\newtheorem{remark}[definition]{Remark}

\newcommand{\refprop}[1]{Proposition~\ref{#1}}

\newcommand{\reflemma}[1]{Lemma~\ref{#1}}

\newcommand{\refeq}[1]{(\ref{#1})}

\newcommand{\refsec}[1]{Section~\ref{#1}}

\newcommand{\prty}[1]{\bar{#1}}
\newcommand{\VPA}{\mathrm{VPA}}
\newcommand{\VA}{\hbar\mathrm{VA}}
\newcommand{\tldfrX}{\widetilde{\frX}}
\newcommand{\tldX}{\widetilde{X}}
\newcommand{\tldfrU}{\widetilde{\frU}}
\newcommand{\tldU}{\widetilde{U}}
\newcommand{\tldcalD}{\widetilde{\calD}}
\newcommand{\tldmu}{\widetilde{\mu}}
\newcommand{\tldvarphi}{\widetilde{\varphi}}
\newcommand{\confwt}{\mathrm{conf\text{-}wt}}
\newcommand{\Swt}{\mathrm{\bbS\text{-}wt}}
\newcommand{\dE}{{}^{\tau} E} 
\newcommand{\Deg}{d} 
\newcommand{\Abar}{\overline{A}} 
\newcommand{\calAbar}{\overline{\calA}} 
\newcommand{\tldP}{\widetilde{P}} 
\newcommand{\tldH}{\widetilde{H}} 
\newcommand{\tldsfD}{\widetilde{\sfD}}

\newcommand{\defeq}{\ensuremath{\underset{\mathrm{def}}{=}}}
\newcommand{\C}{{\mathbb C}}
\newcommand{\R}{{\mathbb R}}
\newcommand{\Q}{{\mathbb Q}}

\newcommand{\Z}{{\mathbb Z}}

\newcommand{\bbS}{{\mathbb S}}
\newcommand{\bbT}{{\mathbb T}}

\newcommand{\bbV}{{\mathbb V}}
\newcommand{\bbW}{{\mathbb W}}

\newcommand{\bfe}{\mathbf{e}}

\newcommand{\calA}{\mathcal{A}}

\newcommand{\calD}{\mathcal{D}}

\newcommand{\calF}{\mathcal{F}}

\newcommand{\calH}{\mathcal{H}}

\newcommand{\calO}{\mathcal{O}}

\newcommand{\calW}{\mathcal{W}}

\newcommand{\frS}{\mathfrak{S}}

\newcommand{\frU}{\mathfrak{U}}

\newcommand{\frX}{\mathfrak{X}}

\newcommand{\frg}{\mathfrak{g}}

\newcommand{\scW}{\mathscr{W}}

\newcommand{\sfD}{\mathsf{D}}

\DeclareMathOperator{\gEnd}{End}
\DeclareMathOperator{\lEnd}{\operatorname{\mathscr{E}\kern-.1pc\mathit{nd}}}
\DeclareMathOperator{\gHom}{Hom}
\DeclareMathOperator{\lHom}{\operatorname{\mathscr{H}\kern-.1pc\mathit{om}}}

\DeclareMathOperator{\Spec}{Spec}

\DeclareMathOperator{\Proj}{Proj}

\DeclareMathOperator{\Ker}{Ker}
\renewcommand{\Im}{\operatorname{Im}}

\DeclareMathOperator{\Gr}{Gr}

\newcommand{\blkbar}{\raisebox{0.5ex}{\rule{2ex}{0.4pt}}}
\newcommand{\qquot}{/\!\!/}

\newcommand{\bfone}{\mathbf{1}}

\newcommand{\isoto}[1][]{\mathop{\xrightarrow[#1]%
{\rule{0pt}{.9ex}%
{\raisebox{-.4ex}[0ex][-.6ex]{$\mspace{3mu}\sim\mspace{3mu}$}}}}}
\newcommand{\scbul}{{\,\raise1pt\hbox{$\scriptscriptstyle\bullet$}\,}}

\begin{document}
\maketitle

\begin{abstract}
 We construct a family of vertex algebras associated with a family of symplectic
 singularity/resolution, called hypertoric varieties. While the hypertoric varieties are
 constructed by a certain Hamiltonian reduction associated with a torus action,
 our vertex algebras are constructed by (semi-infinite) BRST reduction. The construction
 works algebro-geometrically and we construct sheaves of $\hbar$-adic vertex algebras
 over hypertoric varieties which localize the vertex algebras. We show when the vertex
 algebras are vertex operator algebras by giving explicit conformal vectors. We also 
 show that the Zhu algebras of the vertex algebras, associative $\C$-algebras associated
 with non-negatively graded vertex algebras, gives a certain family of 
 filtered quantizations of the coordinate rings of the hypertoric varieties.
\end{abstract}

\section{Introduction}
\label{sec:intro}

Hypertoric varieties are a family of symplectic singularities and their symplectic resolutions.
They are constructed by Hamiltonian reduction of a symplectic vector space by the action of
a torus, and were originally studied as hyperk\"{a}hler manifolds by R.~Bielawski and 
A.~S.~Dancer in \cite{Bielawski-Dancer}. It is well known that a hypertoric variety $X$
has the universal family of $\C^{\times}$-equivariant Poisson deformations $\tldX$ over 
the vector space $\frg^*$ where $\frg^*$ is the dual of the Lie algebra of the torus 
of the Hamiltonian reduction constructing the hypertoric variety $X$ (See 
\cite{Kaledin-Verbitsky}, \cite{Losev-2016}).
By using quantum Hamiltonian reduction, I.~Musson and M.~Van den Bergh in 
\cite{Musson-Van-der-Burgh} constructed a quantization of the hypertoric varieties,
which we call quantized hypertoric algebras or hypertoric enveloping algebras, and
studied its representation theory. This construction admits a certain localization as
discussed in \cite{Bellamy-K} and \cite{BLPW}. That is, we may construct a sheaf of
noncommutative $\C[[\hbar]]$-algebras over the hypertoric variety whose algebra of
global sections can be identified with the quantized hypertoric algebra. Moreover,
the quantum Hamiltonian reduction can be interpreted as a certain BRST reduction
as studied in \cite{BRST}. In \cite{Losev} and \cite{Losev-2016}, I.~Losev studied
the isomorphism classes of filtered quantizations of the coordinate ring $\C[X]$
of the hypertoric variety $X$ and showed that there existed a universal family of
filtered quantizations of $\C[X]$ by using the result of \cite{Bezrukavnikov-Kaledin}.
Each quantized hypertoric algebra is obtained as a fiber of the universal family of
filtered quantizations.

Affine $\calW$-algebras are a family of vertex algebras which generalizes affine 
vertex algebras associated with affine Lie algebras and the Virasoro vertex algebra.
The affine $\calW$-algebras were constructed by quantized Drinfel'd-Sokolov reduction
in \cite{Feigin-Frenkel} and \cite{Frenkel-Kac-Wakimoto}. The construction can be
interpreted as a certain quantization of Hamiltonian reduction of infinite-dimensional
manifolds. Such a quantization of infinite-dimensional Hamiltonian reduction is called
semi-infinite reduction or (semi-infinite) BRST reduction/cohomologies. Properties of
the BRST cohomologies associated with the quantized Drinfel'd-Sokolov reduction,
including the vanishing of higher cohomologies,
were extensively studied by T.~Arakawa in \cite{Arakawa1}, \cite{Arakawa2}. 
In \cite{AKM}, T.~Arakawa, F.~Malikov and the author introduced the BRST reduction 
for sheaves of $\hbar$-adic vertex algebras over Poisson varieties and showed that
the affine $\calW$-algebras at critical level admitted localization as sheaves over
the corresponding Slodowy varieties. The resulting sheaves of $\hbar$-adic vertex
algebras can be understood as quantization of sheaves of vertex Poisson algebras
called jet bundles over the Slodowy varieties.

In this paper, we construct a new family of vertex algebras and study their structure.
Our construction is based on a semi-infinite BRST reduction associated with the
Hamiltonian reduction constructing the hypertoric varieties.
Moreover, our construction also works for sheaves of $\hbar$-adic vertex algebras, and
our vertex algebras admits a certain localization. Namely, we construct a sheaf of
$\hbar$-adic vertex algebras over the universal family of Poisson deformations $\tldX$
of the hypertoric variety $X$ by using the BRST reduction, and then the vertex algebra
of its global sections coincides with our vertex algebra associated with $X$.
As a corollary of the sheaf-theoretic construction, we describe the vertex algebra by
a certain affine local coordinate of $\tldX$, and show that the sheaf of $\hbar$-adic
vertex algebras is locally isomorphic to the tensor product of a $\beta\gamma$-system and
a Heisenberg vertex algebra (\refprop{prop:local-descr}). 
By this isomorphism, we have a free field realization of
our vertex algebra, which is an analog of the Wakimoto realization for affine vertex 
operator algebras (\refprop{prop:Wakimoto}). 
The vertex algebra may or may not be a vertex operator algebra. We
determine when the vertex algebra is a vertex operator algebra by constructing a
conformal vector when it is a vertex algebra (\refprop{prop:conf-vect-def}).

The Zhu algebra of a $\Z_{\ge 0}$-graded vertex algebra is an associative algebra 
introduced by Y.~Zhu in \cite{Zhu} whose representation theory reflects 
fundamental aspects of the representation theory of the original vertex algebra.
We show that the Zhu algebra of our vertex algebra is a certain family of 
filtered quantizations of the coordinate ring $\C[X]$, which include the 
universal family of quantizations
(\refprop{thm:Zhu-alg}). 

We summarize the content of each section. In \refsec{sec:hypertoric-varieties}, we 
summarize the definition and fundamental properties of the hypertoric varieties.
We explicitly construct certain local coordinates which trivialize the Hamiltonian
reduction in \refsec{sec:triv-Hamil-red} and \refsec{sec:sympl-deform}.
In \refsec{sec:chiral-DQ-alg}, vertex algebras, vertex Poisson algebras
and $\hbar$-adic vertex algebras are introduced. In \refsec{sec:semi-infinite-red},
we introduce the main object of this paper, the semi-infinite BRST reduction 
associated with the hypertoric varieties. In \refsec{sec:Clifford}, we review
the Clifford vertex superalgebras, an ingredient of the BRST cohomology. 
In Sections \ref{sec:Poisson-brst}--\ref{sec:zero-Poisson-BRST}, we construct the
jet bundle over a hypertoric variety by the BRST reduction. The results in these
sections are used in the following sections. In \refsec{sec:def-brst}, we 
construct a sheaf of $\hbar$-adic vertex algebra over the hypertoric variety by
the BRST reduction. The cochain complex of the BRST reduction is decomposed naturally
into a double complex. In \refsec{sec:conv-double-cpx}, we show that a spectral
sequence associated with the double complex converges to the BRST cohomology.
In \refsec{sec:local-struct}, we study the local structure of the resulting sheaf
of $\hbar$-adic vertex algebras by using the local coordinates in \refsec{sec:triv-Hamil-red}
and \refsec{sec:sympl-deform}.
In \refsec{sec:F-action}, we construct a vertex algebra from the $\hbar$-adic vertex
algebra of global sections of our sheaf by using a certain symmetry of equivariant torus
action on the sheaf of $\hbar$-adic vertex algebras. We call the obtained vertex algebra
a hypertoric vertex algebra. We also construct an analog of Wakimoto realization in
\refsec{sec:F-action}.
In \refsec{sec:conf-vect}, we determine when the hypertoric vertex algebra is a vertex
operator algebra and construct its conformal vector if it is. 
Finally, in \refsec{sec:Zhu-algebra}, We consider the Zhu algebra of the hypertoric vertex algebra.

\subsection*{Acknowledgments}
A primitive idea for the construction of the hypertoric vertex algebras arose from
discussions with Tomoyuki Arakawa about localization of affine $\calW$-algebras \cite{AKM}.
The author is deeply grateful to Tomoyuki Arakawa for numerous discussions and suggestions.
The author also thanks to Yoshihisa Saito, Naoki Genra, Ryo Sato and Hironori Oya for
valuable comments.

\section{Preliminaries}
\label{sec:preliminaries}

Let $G$ be a torus and $V$ be a $G$-module. We denote the subset of
all $G$-invariant elements of $V$ by $V^G$. For a character 
$\theta : G \longrightarrow \C^*$, we denote the subset of all
$G$-semi-invariant elements of weight $\theta$ by
$V^{G, \theta}$. For a fractional character $\theta \in \gHom(G, \C^{\times}) \otimes_{\Z} \Q$
we also consider the space $V^{G, \theta}$ but it is zero unless $\theta \gHom(G, \C^{\times})$.
For an element $v \in V$, let 
$G_v = \{ g \in G \;\vert\; g \cdot v = v \}$ be the stabilizer of $v$.

For a commutative algebra $A$ over $\C$, let $\Spec A$ be the affine 
scheme associated with $A$. For a commutative graded algebra 
$A = \bigoplus_{n \in \Z_{\ge 0}} A_n$, let $\Proj A$ be the projective 
scheme over $\Spec A_0$, which is associated with $A$. Throughout the paper,
we only consider integral, separated and reduced schemes over $\C$. 
We call them varieties.

Let $X$ be a variety over $\C$.
For a sheaf $\calF$ on $X$ and an open subset $U \subset X$, we denote
the set of local sections of $\calF$ on $U$ by $\calF(U)$ or $\Gamma(U, \calF)$.
We denote the structure sheaf of $X$ by $\calO_X$ and the coordinate ring
of $X$ by $\C[X] = \calO_X(X)$.

\section{Hypertoric varieties}
\label{sec:hypertoric-varieties}

In this section, we recall the definition and fundamental properties of hypertoric
varieties. The definition is given by Hamiltonian reduction by an action of
a torus on a symplectic vector space. We will follow the algebraic presentation given in
\cite{Hasel-Sturmfels}. We consider the same setting as one in \cite{Bellamy-K}, and refer
it for detail of our setting.

\subsection{Hamiltonian torus action}
\label{sec:hamilton-act}

Fix positive integers $1 \le M < N$. Let $V = \C^N$ be an $N$-dimensional vector space,
and let $G = (\C^{\times})^{M}$ be a $M$-dimensional torus. We consider that $G$
acts algebraically on $V$ and take a basis of $V$ such that the corresponding coordinate
functions $x_1$, $\dots$, $x_N \in V^* \subset \C[V]$ are weight vectors with respect to
the action of $G$.
Then, the action of $G$ is given by a $M \times N$ integer-valued matrix 
$\Delta = (\Delta_{ij})_{1\le i \le M, 1 \le j \le N}$ as 
$(t_1, \dots, t_M) \cdot x_j = t_1^{\Delta_{1j}} \dots t_M^{\Delta_{Mj}} x_j$ for 
$(t_1, \dots, t_M) \in G = (\C^{\times})^M$. Setting $\Delta_j = (\Delta_{ij})_{i=1, \dots, M}$,
the $j$-th column of the matrix $\Delta$, $\Delta_j$ is the weight of $x_j$ with respect to 
the $G$-action. 
 We assume that $M \times M$ minors of $\Delta$
are relatively prime. This ensures that the map $\Z^N \xrightarrow{\Delta} \Z^M$ is surjective
and hence the stabilizer of a generic point is trivial.

The action of $G$ on $V$ induces an action on its cotangent bundle $T^*V = V \oplus V^*$.
Set $\Delta^{\pm} = (\Delta, -\Delta)$, a $M \times 2N$ matrix, and let $y_1$, $\dots$, 
$y_N \in V \subset \C[V^*]$ be dual to $x_1$, $\dots$, $x_N$. Then, the action of $G$
on $T^*V$ is given by the matrix $\Delta^{\pm}$ as the action on $V$ is given by $\Delta$.
We consider that $T^*V$ is a symplectic vector space with the standard symplectic form
$\omega = dx_1 \wedge dy_1 + \dots + dx_N \wedge dy_N$. Then, the action of $G$ on $T^*V$
is Hamiltonian and we have a moment map $\mu: T^*V \longrightarrow \frg^*$ given by
\begin{equation}
\label{eq:moment-map}
 \mu((x_1, \dots, x_N, y_1, \dots, y_N)) = \Bigl( \sum_{j=1}^{N} \Delta_{ij} x_j y_j \Bigr)_{1 \le i \le M},
\end{equation}
where $\frg = Lie\,G = \C^M$ is the Lie algebra of $G$. Let $A_1$, $\dots$, $A_M$
be the standard basis of $\frg = \C^M$. The moment map $\mu$ induces a linear map
\[
 \mu^*: \frg = \bigoplus_{i=1}^M \C A_i \longrightarrow \C[T^*V], \qquad 
A_i \mapsto \sum_{j=1}^{N} \Delta_{ij} x_j y_j,
\]
which we call the comoment map. By using the Poisson bracket $\{\blkbar, \blkbar\}$
on the structure sheaf $\calO_{T^* V}$ of the symplectic space $T^* V$, the induced
$\frg$-action on $\calO_{T^* V}$ is described by the comoment map; namely, 
an element $A \in \frg$ acts on $\calO_{T^*V}$ by $A \mapsto \{\mu^*(A), \blkbar\}$.

\subsection{Stability condition}
\label{sec:stability-cond}

We identify $\Q^M$ with the space of fractional characters $\gHom(G, \C^{\times}) \otimes_{\Z} \Q$
of the torus $G = (\C^{\times})^M$. We fix $\delta \in \Q^M$ which we call a stability 
parameter. 

Let $S \subset T^*V$ be a subvariety of $T^*V$ which is closed under the action of $G$.
A point $p \in S$ is called $\delta$-semistable if there exists an $m \in \Z_{> 0}$ such
that we have a function $f \in \C[S]^{G, m\delta}$ with $f(p) \ne 0$. A point $p \in S$
is called $\delta$-stable if, in addition, its stabilizer $G_p$ is finite. We denote the 
subset of all $\delta$-semistable points $S^{ss}_{\delta}$ or simply $S^{ss}$. Also 
the set of all $\delta$-stable points are denoted $S^{s}_{\delta}$ or simply $S^{s}$.
The stability parameter $\delta$ is said to be effective if $S^{ss}_{\delta} \ne \emptyset$.
We say that two effective stability parameters $\delta_1$, $\delta_2$ such that 
$S^{ss}_{\delta_1}$, $S^{ss}_{\delta_2} \ne 0$ are equivalent if
$S^{ss}_{\delta_1} = S^{ss}_{\delta_2}$. In the above situation, we have a rational 
polyhedral fan $\Delta(G, S)$ in $\Q^M$, called the G.I.T. fan, whose
support is the set of all effective parameters $\delta$ such that $S^{ss}_{\delta} \ne 0$
and whose walls are given by all stability parameters $\delta$ such that 
$S^{s}_{\delta} \ne S^{ss}_{\delta}$. Under our assumption on the matrix $\Delta$, the maximal
cones of $\Delta(G, T^*V)$ are $M$-dimensional. We call such cones $M$-cones. 
The matrix $\Delta$ is said to be unimodular if every $M \times M$ minor of $\Delta$ takes values in
$\{-1, 0, 1\}$. 

\subsection{Definition of hypertoric varieties}
\label{sec:def-hypertoric-var}

Now we define hypertoric varieties. 
Fix an effective stability parameter $\delta \in \Q^M$,
and let $\frX = (T^* V)^{ss}_\delta \subset T^* V$ be the subset of all 
$\delta$-semistable points of $T^* V$.
For any $\chi \in \frg^*$, 
the level set $\mu^{-1}(\chi)$ of level $\chi$ 
with respect to the moment map $\mu: T^*V \longrightarrow \frg^*$ is
closed under the action of $G$.  For a subset $S \subset T^*V$, 
two points $p$, $q \in S$ are said to be S-equivalent if the closed $G$-orbits 
$\overline{G \cdot p}$ and $\overline{G \cdot q}$ intersect in $S$.

Then, we define a hypertoric variety
associated with the action of $G$ and the stability parameter $\delta$ as follows:
\begin{definition}
 \label{def:hypertoric-var}
 A hypertoric variety $X_{\delta}$ associated with the action of $G$ on $T^*V$ and 
 the stability parameter $\delta$ is given by the quotient space 
\[
 X_{\delta} = (\mu^{-1}(0) \cap \frX) / \sim
\]
 where $\sim$ is the S-equivalence.
\end{definition}

Recall that the $G$-action on $T^* V$ induces an action of $G$ on the structure sheaf
$\calO_{T^* V}$.
By the fundamental fact of the geometric invariant theory, the hypertoric variety 
$X_{\delta}$ is constructed
as a projective scheme over $X_0 \defeq \Spec \C[\mu^{-1}(0)]^G$;
\begin{equation}
\label{eq:proj-scheme} 
 X_\delta \simeq \Proj \bigoplus_{m \in \Z_{\ge 0}} \C[\mu^{-1}(0)]^{G, m\delta}.
\end{equation}
In the following, we summarize fundamental properties of hypertoric varieties.

\begin{proposition}[\cite{Hasel-Sturmfels}, Proposition 6.2; see also \cite{Bellamy-K}, Corollary~4.13]
 \label{prop:smooth-hypertoric}
 If $\delta$ is in the interior of a $M$-cone of $\Delta(G, \mu^{-1}(0))$ then the
 hypertoric variety $X_{\delta}$ is an orbifold. It is smooth if and only if $\delta$
 is in the interior of a $M$-cone of $\Delta(G, \mu^{-1}(0))$ and $\Delta$ is unimodular.
 Moreover, the walls of the G.I.T. fan are $\sum_{j \in J} \Q \, \Delta_j$ where 
 $J \subset \{1, \dots, N\}$ is any subset such that
 $\dim_{\Q}(\sum_{j \in J} \Q \,\Delta_j) = M-1$.
\end{proposition}

\begin{lemma}[\cite{Bellamy-K}, Lemma 4.7]
\label{lemma:flat-moment} 
 The moment map $\mu$ is flat and $\mu^{-1}(0)$ is a reduced complete intersection in
 $T^*V$.
\end{lemma}

In the rest of the paper, we fix a unimodular matrix $\Delta$ and an effective stability 
parameter $\delta$ which lies in the interior of a $M$-cone of $\Delta(G, \mu^{-1}(0))$.
By \refprop{prop:smooth-hypertoric}, for such $\Delta$ and $\delta$, we have a resolution of
singularity $X_{\delta} \longrightarrow X_{0}$. We also
denote $X_{\delta}$ simply $X$. We denote the morphism of the resolution $\pi$; i.e.
$\pi: X \longrightarrow X_0$. Note that the symplectic structure on $T^*V$ induces
a symplectic structure on $X$. It also induces a Poisson structure on $X_0$ and the
morphism $\pi$ preserves these Poisson structures; i.e. we have a homomorphism of
Poisson algebras $\calO_{X_{0}} \longrightarrow \calO_{X}$.

Now we consider a certain basic fact for the semistable locus $\frX = (T^*V)^{ss}_{\delta}$
with respect to the stability condition $\delta \in \Q^M$. 

In the rest, we identify the space of fractional parameters 
$\gHom(G, \C^{\times}) \otimes_{\Z} \Q$, its dual space and $\Q^M$. We also identify 
the natural pairing between these spaces and the standard inner product of $\Q^M$, and
denote them $( \blkbar, \blkbar )$. We denote the set of common zeros 
of the polynomials $f_1$, $\dots$, $f_r$ $V(f_1, \dots, f_r) \subset T^*V$.

\begin{lemma}[\cite{Bellamy-K}, Lemma~4.3]
 \label{lemma:BK-lemma-4.3}
 For a point $p \in T^*V$, we set the subsets of indices
 $J_1 = \{ j \in \{1, \dots, N\} \;|\; x_j(p) \ne 0 \}$ and 
 $J_2 = \{ j \in \{1, \dots, N\} \;|\; y_j(p) \ne 0 \}$.
 Then, $p \in (T^*V)^{ss}_{\delta}$ if and only if 
 $\delta \in \sum_{j \in J_1} \Q_{\ge 0} \,\Delta_j + \sum_{j \in J_2} \Q_{\le 0} \,\Delta_j$.
\end{lemma}

\begin{proposition}
 \label{prop:semistab-locus}
 For a semistable point $p \in \frX = (T^*V)^{ss}_{\delta}$, there exists a subset of indices
 $\{j_1, \dots, j_M\} \subset \{1, \dots, N\}$ such that $x_{j_i}(p) \ne 0$ or 
 $y_{j_i}(p) \ne 0$ for any $i=1$, $\dots$, $M$ 
 and $\det(\Delta_{j_1}, \dots, \Delta_{j_M}) = \pm 1$.
\begin{proof}
 Set $J_1 = \{ j \;|\; x_j(p) \ne 0\}$, $J_2 = \{ j \;|\; y_j(p) \ne 0\}$ and
 $J = J_1 \cup J_2$.
 By \reflemma{lemma:BK-lemma-4.3}, $p$ is $\delta$-semistable if and only if
 $\delta \in \sum_{j \in J_1} \Q_{\ge 0} \Delta_j + \sum_{j \in J_2} \Q_{\le 0} \Delta_j$.
 Thus, we have $\delta \in \sum_{j \in J} \Q \Delta_j \subset \Q^M$. Since we assume that
 the hypertoric variety $X_{\delta}$ is smooth, we have $\sum_{j \in J} \Q \Delta_j = \Q^M$
 by \refprop{prop:smooth-hypertoric}. 
 Take $j_1, \dots, j_M \in J$ such that $\Delta_{j_1}$, $\dots$, $\Delta_{j_M}$ are linearly
 independent. Then, $\det(\Delta_{j_1}, \dots, \Delta_{j_M}) = \pm 1$ because the matrix
 $\Delta$ is unimodular.
\end{proof}
\end{proposition}

\subsection{Local trivialization of Hamiltonian reduction}
\label{sec:triv-Hamil-red}

Now we construct an affine open covering of $\frX$ which trivializes the Hamiltonian 
reduction with respect to $G$. 

Fix a subset of indices $J = \{j_1, \dots, j_M \}\subset \{1, \dots, N\}$ such that
the minor $\det (\Delta_{j_1}, \dots, \Delta_{j_M}) = \pm 1$. We set 
\[
 \frU_J = \{p \in \frX \;|\; x_j(p) \ne 0 \text{ or } y_j(p) \ne 0 \text{ for any } j \in J \}.
\]
By \refprop{prop:semistab-locus}, we have $\frX = \bigcup_{J} \frU_J$. The stability parameter
$\delta$ can be written in a linear combination of $\{\Delta_{j}\}_{j \in J}$:
$\delta = \sum_{j \in J} \alpha_j \Delta_{j}$ where $\alpha_{j} \in \Q$. Note that 
$\alpha_j \ne 0$ for all $j \in J$ since otherwise $\delta$ lies on the G.I.T. walls 
by \refprop{prop:smooth-hypertoric}. 
Set $J_1 = \{ j \;|\; \alpha_j > 0 \}$ and $J_2 = \{ j \;|\; \alpha_j < 0 \}$. Then,
by \reflemma{lemma:BK-lemma-4.3}, 
we have $x_j(p) \ne 0$ for $j \in J_1$, $y_j(p) \ne 0$ for $j \in J_2$ and $J = J_1 \sqcup J_2$.
Thus we have the following finer description; 
\begin{equation}
 \label{eq:trivialize-open}
  \frU_J = \{ p \in \frX \;|\; x_j(p) \ne 0 \text{ for } j \in J_1,\; y_j(p) \ne 0 \text{ for }
j \in J_2 \}.
\end{equation}

We show that the Hamiltonian reduction with respect to the $G$-action is trivialized 
locally on each open subset $\frU_J$. 
By multiplying a certain positive integer 
$m \in \Z_{> 0}$ to $\delta = \sum_{j \in J} \alpha_j \Delta_{j}$, we have 
$m \delta = \sum_{j \in J} (m \alpha_j) \Delta_{j}$ so that $m \alpha_j \in \Z$. 
Since the weight of $x_j$ with respect to the $G$-action is $\Delta_j$, we have a polynomial
\[
 f_J = \prod_{j \in J_1} x_{j}^{m \alpha_j} \prod_{j \in J_2} y_{j}^{- m \alpha_j}
 \in \C[\frX]^{G, m \delta}
\]
of weight $m \delta$ such that $f_J(p) \ne 0$ for any $p \in \frU_J$. Note that
$x_j^{-1} = f_J^{-1} (f_J x_j^{-1}) \in \calO_{\frX}(\frU_J)$ 
(resp. $y_j^{-1} = f_J^{-1} (f_J y_j^{-1}) \in \calO_{\frX}(\frU_J)$) for 
$j \in J_1$ (resp. $j \in J_2$). Since 
$\det(\Delta_{j})_{j \in J} = \pm 1$,  for each 
$i=1$, $\dots$, $M$, there exist $\lambda_{ij} \in \Z$ for $j \in J$ such that 
$\sum_{j \in J} \lambda_{ij} \Delta_j = \bfe_i$
where $\bfe_i$ is the $i$-th standard basis of $\Z^M$. Set
\[
 T_i^{J} = \prod_{j \in J_1} x_j^{\lambda_{ij}} \prod_{j \in J_2} y_j^{- \lambda_{ij}} \in \calO_{\frX}(\frU_J). 
\]
Then, $T^J_i$ is a local section of weight $\bfe_i$ with respect to the $G$-action and 
it is invertible in $\calO_{\frX}(\frU_J)$ for $i=1$, $\dots$, $M$. In the following,
we also write simply $T_i$ when there is no chance to confuse.
For each $j \not\in J$, we have $G$-invariant local sections
\[
 a_j^{*J} = x_j T_1^{- \Delta_{1j}} \cdots T_M^{- \Delta_{Mj}}, \quad
 a_j^J = y_j T_1^{\Delta_{1j}} \cdots T_M^{\Delta_{Mj}} \in \calO_{\frX}(\frU_J).
\]
Again we also write simply $a_j^*$, $a_j$ instead of $a^{*J}_j$, $a^J_j$ when there is no
confusion.
Note that $\{j\}$ for $j \not\in J$, $J_1$ and $J_2$ are disjoint with one another, 
and hence $T_1$, $\dots$, $T_M$ contain at most one from each symplectic pair $(x_k, y_k)$
for $k=1$, $\dots$, $N$. Thus we have 
$\{ a_j, a_j^* \} = \{ y_j, x_j \} = 1$ for $j \not\in J$ and 
$\{ T_i, a^*_j \} = \{ T_i, a_j \} = \{T_i, T_{i'} \} = 0$ for 
$i$, $i'=1$, $\dots$, $M$ and $j \not\in J$.

Now we describe the trivialization of the Hamiltonian reduction locally on $\frU_J$.
For $i=1$, $\dots$, $M$, put 
$\gamma_i = \mu^*(A_i) = \sum_{j=1}^{N} \Delta_{ij} x_j y_j \in \calO_{\frX}(\frX)$.
Then we have an identity
\begin{equation}
\label{eq:local-triv-Hamil-1} 
 \calO_{\frX}(\frU_J) = 
\C[\; a_j^*, a_j \;|\; j \not\in J\;] \otimes_{\C} \C[\, T_1^{\pm}, \dots, T_M^{\pm}\,]
\otimes_{\C} \C[\,\gamma_1, \dots, \gamma_M\,].
\end{equation}

Indeed, we can describe the generators of $\calO_{\frX}(\frU_J)$ as polynomials of the
generators in the right hand side.
For $j \not\in J$, we have
\[
 x_j = a_j^* T_1^{\Delta_{1j}} \dots T_M^{\Delta_{Mj}}, \quad
 y_j = a_j T^{- \Delta_{1j}} \dots T_M^{- \Delta_{Mj}}.
\]
Note that $a_j^* a_j = x_j y_j$ for $j \not\in J$. Since the matrix 
$(\Delta_{ij})_{i=1, \dots, M, j \in J}$ is invertible with the inverse matrix
$(\lambda_{ji})_{j \in J, i=1, \dots, M}$, from the identity
$\gamma_i - \sum_{j \not\in J} \Delta_{ij} a_j^* a_j = \sum_{j \in J} \Delta_{ij} x_j y_j$,
we obtain
$x_j y_j = \sum_{i=1}^{M} \lambda_{ij} (\gamma_i - \sum_{k \not\in J} \Delta_{ik} a_k^* a_k)$
for $j \in J$.
Thus, for $j \in J_1$, we have
\[
 x_j = T_1^{\Delta_{1j}} \dots T_M^{\Delta_{Mj}}, \qquad
 y_j = T_1^{- \Delta_{1j}} \dots T_M^{- \Delta_{Mj}} \sum_{i=1}^{M} \lambda_{ij}\Bigl(\gamma_i - \sum_{k \not\in J} \Delta_{ik} a_k^* a_k\Bigr).
\]
It is clear that a similar identity holds for $j \in J_2$. This implies that 
the identity \refeq{eq:local-triv-Hamil-1} holds.

Note that we have 
$\{\gamma_i, a_j^*\} = \{\gamma_i, a_j\} = 0$ for $i=1$, $\dots$, $M$, $j \not\in J$,
and $\{\gamma_i, T_j\} = T_j$ for $i$, $j = 1$, $\dots$, $M$ by the construction. 
We regard $\gamma_1$, $\dots$, $\gamma_M$ as a linear basis of the Lie algebra $\frg$ 
through the homomorphism $\mu^*$. Then, the identity \refeq{eq:local-triv-Hamil-1} gives
an isomorphism of Poisson algebras
\begin{equation}
 \label{eq:local-triv-Hamil-2}
 \calO_{\frX}(\frU_J) \simeq
  \C[T^* \C^{N-M}] \otimes_{\C} \C[G] \otimes_{\C} \C[\frg^*] \simeq 
  \C[T^*\C^{N-M}] \otimes_{\C} \C[T^* G]
\end{equation}
and thus we have the trivialization 
$\frU_J \simeq T^* \C^{N-M} \times T^* G$. Set $U_J = (\mu^{-1}(0) \cap \frU_J) / G \subset X$.
Since the $G$-action and the moment map $\mu$ are trivialized, we have 
$U_J \simeq T^* \C^{N-M}$ as symplectic manifolds. Then, we have 
an affine open covering $X = \bigcup_{J} U_J$ with Darboux coordinate 
$(a^{*J}_j, a^J_j)_{j \not\in J}$ for each $J$. 

We denote the trivialization $\nu_J: \frU_J \longrightarrow U_J \times T^* G$
and the corresponding isomorphism
$\nu^*_J: \calO_{U_J} \otimes \calO_{T^* G} \longrightarrow \calO_{\frX}\vert_{\frU_J}$. 
For $I$ and $J$, we denote the coordinate transformation 
$\varphi_{IJ} = \nu_I \circ \nu_J^{-1}: U_J \times T^* G \longrightarrow U_I \times T^* G$, and
the corresponding isomorphism
$\varphi^*_{IJ} = (\nu^*_J)^{-1} \circ \nu^*_I: \calO_{U_I} \otimes \calO_{T^* G} \longrightarrow \calO_{U_J} \otimes \calO_{T^* G}$ on $U_I \cap U_J$.
This induces the coordinate transformation 
$\varphi_{IJ} : U_J \times G \longrightarrow U_I \times G$ and the corresponding isomorphism
$\varphi^*_{IJ} : \calO_{U_I} \otimes \calO_G \longrightarrow \calO_{U_J} \otimes \calO_G$
because $\gamma_1$, $\dots$, $\gamma_M$ are global sections.
Note that, the isomorphism $\varphi^*_{IJ} : \calO_{U_I} \longrightarrow \calO_{U_J}$ coincides
with the coordinate transformation between $\calO_X \vert_{U_I} = \calO_{U_I}$ and 
$\calO_X \vert_{U_J} = \calO_{U_J}$ since it is the coordinate translation of the $G$-torsor
$\mu^{-1}(0) \cap \frX \longrightarrow X$.

\subsection{Symplectic deformation of the hypertoric variety $X$}
\label{sec:sympl-deform}

For the symplectic variety $X$, it is known that there exists
a universal family of filtered Poisson deformations of the symplectic structure of $X$,
which explicitly given as follows.

Set $\widetilde{\frX} = \frX \times \frg^*$. 
We regard $\tldfrX$ as a smooth algebraic Poisson variety where $\frg^*$ is equipped
with the trivial Poisson structure.
We extend the moment map $\mu$ to 
$\widetilde{\mu} : \tldfrX \longrightarrow \frg^*$ such that the corresponding 
comoment map 
$\widetilde{\mu}^*: \frg \longrightarrow \C[\widetilde{\frX}] = \C[\frX] \otimes \C[\frg^*]$ 
is given by $\widetilde{\mu}^*(A_i) = \mu^*(A_i) - c_i$ where we denote the standard basis 
of $\frg \subset \C[\frg^*]$ by $c_1$, $\dots$, $c_M$ instead of $A_1$, $\dots$, $A_M$ in
order to avoid confusion. Clearly the torus $G$ acts on $\widetilde{\frX}$ freely and 
the $G$-action preserves the preimage $\widetilde{\mu}^{-1}(0)$. Then, we define 
the Poisson manifold
$\widetilde{X} = \widetilde{X}_{\delta} = \widetilde{\mu}^{-1}(0) / G \simeq \frX / G$.
Here the last isomorphism is induced from the obvious isomorphism 
$\widetilde{\mu}^{-1}(0) \simeq \frX$ which identifies $c_i$ with $\mu^{*}(A_i)$
for $i=1$, $\dots$, $M$. By the second projection 
$\rho: \widetilde{\frX} = \frX \times \frg^* \longrightarrow \frg^*$ induces the morphism
$\rho: \widetilde{X} \longrightarrow \frg^*$ of Poisson schemes, and we have
$\rho^{-1}(0) \simeq X$. Note that $\widetilde{X}$ is a symplectic scheme over $\frg^*$ 
and the isomorphism $\rho^{-1}(0) \simeq X$ is an isomorphism of holomorphic symplectic
manifold. It is known that $\tldX$ is a universal family of filtered Poisson deformations
of $X$ over $\frg^* \simeq H^2(X, \C)$, namely, the structure sheaf $\calO_{\tldX}$
is a universal family of filtered Poisson deformations of the sheaf of Poisson algebras 
$\calO_X$. Moreover, the family is equivariant with respect to an action of a torus 
$\bbS = \C^{\times}$ which we discuss in \refsec{sec:F-action}. Refer \cite{Losev}
for the universality of the above $\C^{\times}$-equivariant Poisson deformations,
which is based on results of \cite{Kaledin-Verbitsky}.

While the hypertoric varieties $X$ and $\tldX$ are constructed by Hamiltonian reduction by
the action of the torus $G$, their structure sheaves can be constructed also by Hamiltonian
reduction of algebras.
Namely, The structure
sheaf of $\tldX$ is given by the following (dual) Hamiltonian reduction
\[
 \calO_{\widetilde{X}} \simeq 
\Bigl(p_{*}\bigl(\calO_{\tldfrX} \Bigm/ \sum_{i=1}^{M} 
\calO_{\tldfrX} \tldmu^*(A_i) \bigr)
\Bigr)^{G}
=
\Bigl(p_{*}\bigl(\calO_{\tldfrX} \Bigm/ \sum_{i=1}^{M} 
\calO_{\tldfrX} (\mu^*(A_i) - c_i) \bigr)
\Bigr)^{G}
\]
where $p: \tldfrX \longrightarrow \tldX$ is the projection.
It is an algebra over $\C[c_1, \dots, c_M] = \C[\frg^*]$. 
The hypertoric variety $X$
is the fiber of $\tldX \longrightarrow \frg^*$ at 
$c_1 = \dots = c_M = 0$, and we have
\[
 \calO_{X} = 
\Bigl(p_{*}\bigl(\calO_{\frX} \Bigm/ \sum_{i=1}^{M} 
\calO_{\frX} \mu^*(A_i) \bigr)
\Bigr)^{G}.
\]

We consider local trivialization of the Hamiltonian reduction of $\widetilde{\frX}$ by 
the $G$-action. Recall the affine open covering $\frX = \bigcup_J \frU_J$ which trivializes
the Hamiltonian reduction in \refsec{sec:triv-Hamil-red}. Set
$\widetilde{\frU}_J = \frU_J \times \frg^* \subset \tldfrX$ for each $J$. Then, we have 
an affine open covering $\widetilde{\frX} = \bigcup_{J} \widetilde{\frU}_J$. 
Since the $G$-action preserves $\frU_J$ and it acts on $\frg^*$ trivially, $\widetilde{\frU}_J$
is also preserved by the $G$-action. We set 
$\widetilde{U}_J = \widetilde{\mu}^{-1}(0) \cap \widetilde{\frU}_J / G$, and then we have
an open covering $\widetilde{X} = \bigcup_{J} \widetilde{U}_J$.

By the 
trivialization of the Hamiltonian reduction on $\frU_J$ discussed in \refsec{sec:triv-Hamil-red},
we have an isomorphism 
$\widetilde{\frU}_J \simeq T^* \C^{N-M} \times G \times \frg^* \times \frg^*$.
The isomorphism is given by the following description of the algebra of local sections 
$\calO_{\widetilde{\frX}}(\widetilde{\frU}_J)$:
\begin{equation}
\label{eq:tld-local-trivial} 
 \calO_{\widetilde{\frX}}(\widetilde{\frU}_J) \simeq 
\C[\; a_j^*, a_j \;|\; j \not\in J\;] \otimes_{\C} \C[\, T_1^{\pm}, \dots, T_M^{\pm}\,]
\otimes_{\C} \C[\,\gamma_1, \dots, \gamma_M\,] 
\otimes_{\C} \C[\, c_1, \dots, c_M \,]
\end{equation}
where the local sections $a^*_j$, $a_j$, $T_i$, $\gamma_i$ are defined 
in \refsec{sec:triv-Hamil-red}. In the above local coordinate, the comoment map 
$\widetilde{\mu}^*$ is given by $\widetilde{\mu}^*(A_i) = \gamma_i - c_i$ for $i=1$, $\dots$,
$M$. Moreover, since the $G$-action on $\calO_{\widetilde{\frX}}(\widetilde{\frU}_J)$ 
corresponds to
the $\frg$-action $A_i \mapsto \{ \widetilde{\mu}^*(A_i), \blkbar \}$, the torus $G$ acts on 
$\C[a^*_j, a_j \;\vert\; j \not\in J]$, $\C[\gamma_1, \dots, \gamma_M]$ and $\C[c_1, \dots, c_M]$
trivially, and $T_i$ has weight $\bfe_i$ with respect to the $G$-action for $i=1$, $\dots$, $M$.
Therefore, we have 
\begin{equation}
 \label{eq:tld-triv-G-torsor}
  \calO_{\tldmu^{-1}(0)}(\tldfrU_J) \simeq  \C[\,a^*_j, a_j \;\vert\; j \not\in J \,] 
  \otimes_{\C} \C[\, T_1^{\pm}, \dots, T_M^{\pm} \,]\otimes_{\C} \C[\,c_1, \dots, c_M\,].
\end{equation}
and
\begin{equation}
\label{eq:tld-local-coord} 
\calO_{\widetilde{X}}(\widetilde{U}_J) = 
 \calO_{\widetilde{\mu}^{-1}(0)}(\widetilde{\frU}_J)^G \simeq 
 \C[\,a^*_j, a_j \;\vert\; j \not\in J \,] \otimes_{\C} \C[\,c_1, \dots, c_M\,].
\end{equation}

It induces the isomorphism $\widetilde{U}_J \simeq T^*\C^{N-M} \times \C^M$, and hence 
the open covering $\widetilde{X} = \bigcup_J \widetilde{U}_J$ is an affine open covering.
Let $\widetilde{\nu}_J : \widetilde{\frU}_J \longrightarrow \widetilde{U}_J \times G \times \frg^*$
be the above trivialization, and we denote the corresponding algebra isomorphism 
$\widetilde{\nu}^*_J: \calO_{\widetilde{U}_J} \otimes_{\C} \calO_{G} \otimes_{\C} \calO_{\frg^*} \longrightarrow \calO_{\widetilde{\frU}_J}$.
Then we have the coordinate transformation over $\tldfrU_I \cap \tldfrU_J$ for $I$, $J$,
$\widetilde{\varphi}_{IJ} = \widetilde{\nu}_I \circ \widetilde{\nu}_J^{-1}: \widetilde{U}_J \times G \times \frg^* \longrightarrow \widetilde{U}_I \times G \times \frg^*$, 
and the algebra isomorphism 
$\widetilde{\varphi}^*_{IJ} : \calO_{\widetilde{U}_J} \otimes_{\C} \calO_{G} \otimes_{\C} \calO_{\frg^*} \longrightarrow \calO_{\widetilde{U}_I} \otimes_{\C} \calO_{G} \otimes_{\C} \calO_{\frg^*}$.
This coordinate transformation induces the coordinate translation of $G$-torsor 
$\widetilde{\varphi}_{IJ} : \widetilde{U}_J \times G \longrightarrow \widetilde{U}_I \times G$
over $\widetilde{U}_I \cap \widetilde{U}_J$,
and the coordinate translation of local coordinates of $\widetilde{X}$, 
$\widetilde{\varphi}_{IJ} : \widetilde{U}_J \longrightarrow \widetilde{U}_I$. The corresponding
algebra isomorphisms are also denoted $\widetilde{\varphi}^*_{IJ}$.

\section{Sheaves of $\hbar$-adic vertex algebras}
\label{sec:chiral-DQ-alg}

In this section, we review the definitions of vertex algebras and $\hbar$-adic vertex
algebras, and we introduce certain sheaves of vertex Poisson algebras and certain
sheaves of $\hbar$-adic vertex algebras. Based on these sheaves,
we will construct a sheaf of vertex Poisson algebras and a sheaf of $\hbar$-adic
vertex algebras in the next section.

\subsection{Vertex algebras and $\hbar$-adic vertex algebras}
\label{sec:vertex-algebras}

A vertex algebra $V$ is a vector space over $\C$ equipped with the following
structure; the vacuum vector $\bfone \in V$, the translation operator $\partial: V \longrightarrow V$
and the vertex operator 
$Y(a, z) = \sum_{n \in \Z} a_{(-n-1)} z^{n} \in \gEnd_{\C}(V)[[z, z^{-1}]]$ 
for each $a \in V$ subject to the following axioms:
\begin{enumerate}
 \item $Y(a, z)$ is linear with respect to $a \in V$,
 \item $Y(a, z)$ is a field, i.e. $a_{(n)} b = 0$ for any $a$, $b \in V$ if $n \gg 0$.
 \item $Y(\bfone, z) = Id_V$,
 \item $Y(a, z) \bfone \in V[[z]]$ and $Y(a, z) \bfone \vert_{z=0} = a$ for any $a \in V$,
 \item $[\partial, Y(a, z)] = \partial_z Y(a, z)$ for any $a \in V$, and $\partial \bfone = 0$,
 \item for any $a$, $b \in V$, the vertex operators $Y(a, z)$ and $Y(b, w)$ are mutually
       local; namely, there exists $N \in \Z_{\ge 0}$ such that
       \[
	(z-w)^N [Y(a, z), Y(b, w)] = 0.
       \]
\end{enumerate}
It is well-known that fundamental identities for vertex algebras 
such as $\partial a = a_{(-2)} \bfone$,
$Y(\partial a, z) = \partial_z Y(a, z)$ and the operator product expansion (or so called
Borcherds' identity) follow from the above axioms. We say that the vertex algebra $V$
is commutative if $a_{(n)} = 0$ on $V$ for any $a$ and $n \in \Z_{\ge 0}$.

A vertex Poisson algebra $V$ is a tuple 
$(V, \bfone, \partial, Y_{-}(\blkbar, z), Y_{+}(\blkbar, z))$ where $Y_{-}(\blkbar, z)$,
$Y_{+}(\blkbar, z): V \longrightarrow \gEnd_{\C}(V)[[z, z^{-1}]]$ are fields on $V$, 
\[
 Y_{-}(a, z) = \sum_{n \in \Z_{\ge 0}} a_{(-n-1)} z^n, \quad
 Y_{+}(a, z) = \sum_{n \in \Z_{< 0}} a_{(-n-1)} z^n
\]
such that $(V, \bfone, \partial, Y_{-}(\blkbar, z))$ is a commutative vertex algebra, and
$(V, \partial, Y_{+}(\blkbar, z))$ is a vertex Lie algebra; namely the operators
$a_{(n)}$ satisfy the following relations:
\begin{enumerate}
 \item $a_{(n)} b = (-1)^{n+1} \sum_{j \ge 0} (-1)^j \partial^j (b_{(n+j)} a) / j!$,
 \item $a_{(m)} b_{(n)} c - b_{(n)} a_{(m)} c = \sum_{j \ge 0} \binom{m}{j} (a_{(j)} b)_{(m+n-j)} c$,
 \item $[\partial, Y_{+}(a, z)] = \partial_z Y_{+}(a, z)$, and
 \item $a_{(n)}$ is a derivation with respect to the product ${}_{(-1)}$,
\end{enumerate}
for any $a$, $b$, $c \in V$ and $m$, $n \in \Z_{\ge 0}$.

Let $\hbar$ be an indeterminate, which commutes with any other operators.
An $\hbar$-adic vertex algebra $V$ is a tuple $(V, \bfone, \partial, Y(\blkbar, z))$ such that
$V$ is a flat $\C[[\hbar]]$-module complete in $\hbar$-adic topology, the vacuum vector
$\bfone \in V$ and $\C[[\hbar]]$-linear map $\partial: V \longrightarrow V$ satisfy the same
axiom with the above, and $Y(\blkbar, z): V \longrightarrow \gEnd_{\C}(V)[[z, z^{-1}]]$
is $\C[[\hbar]]$-linear map such that the products ${}_{(n)}$ are continuous with 
respect to $\hbar$-adic topology, and
$(V / \hbar^N V, \bfone, \partial, Y(\blkbar, z))$
is a vertex algebra for each $N \in \Z_{\ge 1}$. Note that a $\hbar$-adic vertex algebra
is not a vertex algebra over $\C[[\hbar]]$ since $Y(a, z)$ is not a field on $V$. Namely
for any $N \in \Z_{\ge 1}$,
$Y(a, z) = \sum_{n \in \Z} a_{(n)} z^{-n-1}$ satisfies $a_{(n)} b \equiv 0$ modulo $\hbar^N$
if $n \gg 0$, but not $a_{(n)} b = 0$. 

Let $(V, \bfone, \partial, Y(\blkbar, z))$ be an $\hbar$-adic vertex algebra. Assume that
$V / \hbar V$ is commutative. Then, $Y_{+}(\blkbar, z) := \hbar^{-1} Y(\blkbar, z)$
modulo $\hbar$ satisfies the axiom of vertex Lie algebras. Thus, 
$(V / \hbar V, \bfone, \partial, Y(\blkbar, z) \bmod \hbar, \hbar^{-1} Y(\blkbar, z) \bmod \hbar)$
is a vertex Poisson algebra.

\subsection{Jet bundles}
\label{sec:jet-bundles}

Let $X$ be a scheme over $\C$. Let $J_{\infty} X$ be the corresponding $\infty$-jet scheme;
i.e. $J_{\infty} X$ is a scheme defined by 
$\gHom(\Spec R, J_{\infty} X) = \gHom(\Spec R[[t]], X)$ for any $\C$-algebra $R$. 
A point of $J_{\infty} X$ represents an $\infty$-jet $p(t) = \sum_{n=0}^{\infty} p_n t^n$
($p_n \in X$) on $X$. A canonical morphism $\pi_{\infty}: J_{\infty} X \longrightarrow X$
is given by $p(t) \mapsto p(0) = p_0$. We consider the direct image of the structure sheaf
$\calO_{J_{\infty} X}$ of the $\infty$-jet scheme $J_{\infty} X$ by the morphism $\pi_{\infty}$.
The obtained sheaf on $X$ is denoted $\calO_{J_{\infty} X}$ by abuse of notation, and
call it the jet bundle on $X$.
The corresponding homomorphism between their structure
sheaves $\pi_{\infty}^*: \calO_X \hookrightarrow \calO_{J_{\infty} X}$ is an injective 
homomorphism of commutative algebras.
The derivation with respect to $t$ on $R[[t]]$ induces a derivation $\partial$ on the jet bundle 
$\calO_{J_{\infty} X}$.
Thus, the jet bundle $\calO_{J_{\infty} X}$ is a sheaf of commutative vertex algebras on $X$.
Moreover, when $X$ is a Poisson scheme, the Poisson bracket $\{\cdot, \cdot\}$ on
$\calO_X$ induces a structure of vertex Poisson algebras on $\calO_{J_{\infty} X}$
satisfying $f_{(0)} g = \{f, g\}$ and $f_{(n)} g = 0$ for $f$, 
$g \in \calO_X \subset \calO_{J_{\infty} X}$ and $n \in \Z_{\ge 1}$. For detail of the 
construction of vertex Poisson algebra structure, see \cite[Lemma 2.1.3.1]{AKM}.

In the present paper, we consider a smooth symplectic manifold $X$. Assume that a local 
Darboux coordinate $(U; x_1, \dots, x_r, y_1, \dots, y_r)$ is given. Then, the algebra of
local sections of the structure sheaf $\calO_X(U)$ is the polynomial ring 
$\C[x_1, \dots, x_r, y_1, \dots, y_r]$ and the Poisson bracket is given by 
$\{y_i, x_j\} = \delta_{ij}$, ${x_i, x_j} = {y_i, y_j} = 0$. The jet bundle looks like 
\[
 \calO_{J_{\infty} X}(U) = \C[x_{1, (-n)}, \dots, x_{r, (-n)}, y_{1, (-n)}, \dots, y_{r, (-n)}
\;|\; n \in \Z_{\ge 1}],
\]
so that we identify $x_{i} = x_{i, (-1)}$, $y_i = y_{i, (-1)}$ under the embedding 
$\calO_{X} \hookrightarrow \calO_{J_{\infty} X}$. The derivation $\partial$ on $\calO_{J_{\infty} X}$
is given by $\partial(a_{(-n)}) = n a_{(-n-1)}$ for $a = x_i$, $y_i$ ($i=1$, $\dots$, $r$) and 
$n \in \Z_{\ge 1}$. Finally, the vertex Poisson algebra structure on $\calO_{J_{\infty} X}(U)$
is given by
\[
Y_{+}(x_{i, (-1)}, z) = - \sum_{n \ge 1} \frac{\partial}{\partial y_{i, (-n)}} z^{-n},
\quad
Y_{+}(y_{i, (-1)}, z) = \sum_{n \ge 1} \frac{\partial}{\partial x_{i, (-n)}} z^{-n}.
\]

\subsection{$\hbar$-adic $\beta\gamma$-systems and $\hbar$-adic Heisenberg vertex algebras}
\label{sec:h-adic-betagamma}

Let $x_1$, $\dots$, $x_N$, $y_1$, $\dots$, $y_N$ be the standard coordinate functions
on $T^*\C^N = \C^{2N}$. We consider that they are Darboux coordinates with respect to
the standard symplectic form. The $\hbar$-adic $\beta\gamma$-system on $\C^{2N} = T^* \C^N$ is
an $\hbar$-adic vertex algebra $\calD^{ch}(\C^{2N})_{\hbar}$ such that 
$\calD^{ch}(\C^{2N})_{\hbar}$ is isomorphic 
\[
 \calD^{ch}(\C^{2N})_{\hbar} = \C[[\hbar]][x_{1,(-n)}, \dots, x_{N, (-n)}, y_{1, (-n)}, \dots,
y_{N, (-n)} \;\vert\; n \in \Z_{\ge 1}] \bfone
\]
as a $\C[[\hbar]]$-module, and its OPEs are given by $x_i(z) y_j(w) \sim - \hbar/(z-w)$,
and $x_i(z) x_j(w) \sim y_i(z) y_j(w) \sim 0$ for $i$, $j=1$, $\dots$, $N$, where
we denote $x_i(z) = Y(x_{i, (-1)}\bfone, z)$ and $y_i(z) = Y(y_{i, (-1)} \bfone, z)$.
Clearly it is an $\hbar$-adic analogue of the vertex algebra $\beta\gamma$-system.

In \cite{AKM}, we discussed localization of algebras of chiral differential operators
(CDOs), including the $\beta\gamma$-system, as sheaves of $\hbar$-adic vertex algebras
on cotangent bundles; i.e. the above $\hbar$-adic $\beta\gamma$-system gives a sheaf
of $\hbar$-adic vertex algebras on $\C^{2N} = T^*\C^N$ as follows: For the $\hbar$-adic 
$\beta\gamma$-system, OPEs (and hence ${}_{(n)}$-products) between vertex operators
are determined by the Wick formula and thus they turn out to be bi-differential operators
in the variables $x_{i, (-n)}$, $y_{i, (-n)}$. Therefore, even for rational functions in
$x_{i, (-n)}$, $y_{i, (-n)}$, the same bi-differential operators give well-defined OPEs 
(${}_{(n)}$-products) between them. Therefore, we have a sheaf of $\hbar$-adic vertex
algebras $\calD^{ch}_{T^* \C^N, \hbar}$ on $T^* \C^N$. See Lemma~2.2.8.1 and 
Theorem~2.2.10.1 in \cite{AKM} for the detail of the above discussion. 

As we discussed in the previous section, the jet bundle $\calO_{J_{\infty} T^* \C^N}$ on 
the symplectic vector space $T^* \C^N$ is equipped with the vertex Poisson algebra 
structure. The $\hbar$-adic $\beta\gamma$-system $\calD^{ch}_{T^*\C^N, \hbar}$ is a
quantization of $\calO_{J_{\infty} T^* \C^N}$; namely, the quotient 
$\calD^{ch}_{T^*\C^N, \hbar} / \hbar \calD^{ch}_{T^*\C^N, \hbar}$ is isomorphic to
$\calO_{J_{\infty} T^* \C^N}$ as vertex Poisson algebras.

Similarly we define an $\hbar$-adic Heisenberg vertex algebra. Let 
$W = \bigoplus_{i=1}^M \C c_i$ be a vector space with a symmetric inner product
$\langle \blkbar, \blkbar \rangle$. Consider the $\hbar$-adic vertex algebra which 
is defined as $\C[[\hbar]]$-module 
\[
 V_{\langle, \rangle, \hbar}(W) = \C[[\hbar]][\,c_{1, (-n)}, \dots, c_{M, (-n)} \;|\; n \in \Z_{\ge 1}\,],
\]
and OPEs are given by $c_i(z) c_j(w) \sim \hbar^2 \langle c_i, c_j \rangle / (z-w)^2$ for 
$i$, $j = 1$, $\dots$, $M$. Clearly, it is a natural $\hbar$-adic analogue of the usual
Heisenberg vertex algebra defined by $(W, \langle\blkbar, \blkbar\rangle)$.
This implies that the Wick formula holds for the OPEs between vertex operators of 
$V_{\langle, \rangle, \hbar}(W)$ and hence the OPEs are defined as bi-differential
operators in the variables $c_{i, (-n)}$ for $i=1$, $\dots$, $M$, $n \in \Z_{\ge 1}$.
Thus, by the same argument for the $\beta\gamma$-system, the $\hbar$-adic vertex algebra
induces a sheaf of $\hbar$-adic vertex algebras on the vector space $W$. We denote the
sheaf $V_{W, \langle, \rangle, \hbar}$. 

\section{Semi-infinite BRST reduction}
\label{sec:semi-infinite-red}

Now we construct a sheaf of $\hbar$-adic vertex algebras on the hypertoric variety 
$\tldX$ in this section. Our construction is based on a vertex algebra analog of 
the Hamiltonian reduction, which we call (semi-infinite) BRST reduction or BRST cohomology.

In \refsec{sec:Clifford}, we introduce an $\hbar$-adic variant of a fermionic vertex 
superalgebra called the Clifford vertex superalgebra or the free field of colored fermions.
To establish fundamental properties of the BRST reduction, we first need to consider 
the corresponding reduction for a sheaf of vertex Poisson algebra, the jet bundle on $\tldX$.
In Sections \ref{sec:Poisson-brst}--\ref{sec:zero-Poisson-BRST}, we introduce the 
BRST reduction for vertex Poisson algebras and study its structure. 
The BRST reduction for a sheaf of $\hbar$-adic vertex algebras is defined in 
\refsec{sec:def-brst}, and we show that the structure of such a sheaf of $\hbar$-adic 
vertex algebras can be studied by using a certain double complex in 
\refsec{sec:conv-double-cpx}.

\subsection{Clifford $\hbar$-adic vertex superalgebra}
\label{sec:Clifford}

In this subsection, we introduce the Clifford $\hbar$-adic vertex superalgebra 
$Cl^{vert}(\frg \oplus \frg^*)$ associated with the vector space $\frg \oplus \frg^*$
with the standard inner product $\langle \blkbar, \blkbar \rangle$. 

We fix a basis $\frg = \bigoplus_{i=1}^{M} \C A_i$ and its dual basis 
$\frg^* = \bigoplus_{i=1}^{M} \C A^*_i$ with respect to $\langle \blkbar, \blkbar \rangle$
as previous sections. Let $\Pi \frg$ (resp. $\Pi \frg^*$) be the odd vector space corresponding
to the even vector space $\frg$ (resp. $\frg^*$), and let 
$\Pi \frg = \bigoplus_{i=1}^{M} \C \psi_i$  (resp. $\Pi \frg^* = \bigoplus_{i=1}^{M} \C \psi^*_i$) 
be the odd basis corresponding to the even basis $\frg = \bigoplus_{i=1}^{M} \C A_i$
(resp. $\frg^* = \bigoplus_{i=1}^{M} \C A^*_i$). We identify the coordinate rings
$\C[\Pi \frg] = \Lambda(\frg^*)$, $\C[\Pi \frg^*] = \Lambda(\frg)$ and 
$\C[T^* \Pi \frg] = \Lambda(\frg \oplus \frg^*)$ where $\Lambda(W)$ is the exterior algebra 
of a vector space $W$. Note that the inner product $\langle \blkbar, \blkbar \rangle$ on
$\frg \oplus \frg^*$ gives a Poisson superalgebra structure on $\Lambda(\frg \oplus \frg^*)$;
$\{\psi_i, \psi^*_j\} = \delta_{ij}$, $\{\psi_i, \psi_j\} = 0 = \{\psi^*_i, \psi^*_j\}$ for
$i$, $j=1$, $\dots$, $M$. A vertex Poisson superalgebra analogue of $\Lambda(\frg \oplus \frg^*)$
(``the jet bundle'' over the super-manifold $T^* \Pi \frg$) is naturally constructed as follows:
Define $\Lambda^{vert}(\frg \oplus \frg^*)$ as an anti-commutative algebra
\[
 \Lambda^{vert}(\frg \oplus \frg^*) = 
\bigwedge_{\substack{i=1, \dots, M \\ n \ge 1}} \C \psi^*_{i (-n)} \wedge
\bigwedge_{\substack{i=1, \dots, M \\ n \ge 1}} \C \psi_{i (-n)} \bfone,
\]
and the Poisson structure is defined by
$\psi_{i (m)}(\psi^*_{j (-n-1)}) = \delta_{m-n, 0} \delta_{ij}$, 
$\psi^*_{i (m)}(\psi^*_{j (-n-1)}) = \delta_{m-n, 0} \delta_{ij}$ and
$\psi_{i (m)}(\psi_{j (-n-1)}) = 0 = \psi^*_{i (m)}(\psi_{j (-n-1)})$ for 
$i$, $j=1$, $\dots$, $M$ and $n$, $m \in \Z_{\ge 0}$. Then $\Lambda^{vert}(\frg \oplus \frg^*)$
is a vertex Poisson superalgebra. Identifying $\psi_i = \psi_{i (-1)} \bfone$, 
$\psi^*_i = \psi^*_{i (-1)} \bfone$ for $i=1$, $\dots$, $M$, the exterior algebra 
$\Lambda(\frg \oplus \frg^*)$ is a subalgebra of $\Lambda^{vert}(\frg \oplus \frg^*)$.

Now we consider the Clifford $\hbar$-adic vertex superalgebra, a quantization of 
the vertex Poisson superalgebra
$\Lambda^{vert}(\frg \oplus \frg^*)$.
Define the $\hbar$-adic vertex superalgebra 
$Cl^{vert}(\frg \oplus \frg^*)$ as a $\C[[\hbar]]$-module,
\[
 Cl^{vert}(\frg \oplus \frg^*) = 
\bigwedge_{\substack{i=1, \dots, M \\ n \ge 1}} \C[[\hbar]] \psi^*_{i (-n)} \widehat{\otimes}
\bigwedge_{\substack{i=1, \dots, M \\ n \ge 1}} \C[[\hbar]] \psi_{i (-n)} \bfone
\]
where $\widehat{\otimes}$ is the completion of the tensor product with respect to the 
$\hbar$-adic topology.
We denote the vertex operators $\psi_i(z) = Y(\psi_i, z)$ and $\psi^*_i(z) = Y(\psi^*_i, z)$.
Then the defining OPEs are given by
\[
 \psi_i(z) \psi^*_j(w) \sim \frac{\hbar \delta_{ij}}{z-w}, \quad
 \psi_i(z) \psi_j(w) \sim 0 \sim \psi^*_i(z) \psi^*_j(w)
\]
for $i$, $j = 1$, $\dots$, $M$. These OPEs give the structure of $\hbar$-adic vertex algebra
on $Cl^{vert}(\frg \oplus \frg^*)$, which we call the Clifford ($\hbar$-adic) vertex superalgebra. 
Clearly we have 
$Cl^{vert}(\frg \oplus \frg^*) / (\hbar) \simeq \Lambda^{vert}(\frg \oplus \frg^*)$ and
thus the Clifford vertex superalgebra $Cl^{vert}(\frg \oplus \frg^*)$
is a quantization of $\Lambda^{vert}(\frg \oplus \frg^*)$.

Note that the vertex Poisson superalgebra $\Lambda^{vert}(\frg \oplus \frg^*)$ and the 
Clifford vertex algebra $Cl^{vert}(\frg \oplus \frg^*)$ are $\Z$-graded by the degree
$\deg(\psi_{i (-n)}) = -1$, $\deg(\psi^*_{i (-n)}) = 1$ and $\deg(\bfone) = 0$ for
$i=1$, $\dots$, $M$ and $n \in \Z$. Let $\Lambda^{vert, n}(\frg \oplus \frg^*)$ 
and $Cl^{vert, n}(\frg \oplus \frg^*)$ be the homogeneous subspaces of degree $n$.
Moreover we have the following decomposition of $\Lambda^{vert, n}(\frg \oplus \frg^*)$
(resp. $Cl^{vert, n}(\frg \oplus \frg^*)$) as a $\C$-vector space (resp. a $\C[[\hbar]]$-module)
\begin{align*}
 \Lambda^{vert, n}(\frg \oplus \frg^*) &= \bigoplus_{p+q=n} \Lambda_{\C}^{vert, p}(\frg^*)
 \otimes \Lambda_{\C}^{vert, q}(\frg), \\
 Cl^{vert, n}(\frg \oplus \frg^*) &= \sideset{}{^{\wedge}}{\bigoplus}_{p+q=n} 
\Lambda_{\C[[\hbar]]}^{vert, p}(\frg^*) \,\widehat{\otimes}\, \Lambda_{\C[[\hbar]]}^{vert, q}(\frg),
\end{align*}
where $\sideset{}{^{\wedge}}{\bigoplus}$ (resp. $\,\widehat{\otimes}\,$) is the completion of 
the direct sum (resp. the tensor product) with respect to the $\hbar$-adic topology, and
\[
 \Lambda_{R}^{vert}(\frg) = 
\bigwedge_{\substack{i=1, \dots, M \\ n \ge 1}} R \psi_{i (-n)}, \quad
 \Lambda_{R}^{vert}(\frg^*) = 
\bigwedge_{\substack{i=1, \dots, M \\ n \ge 1}} R \psi^*_{i (-n)}
\]
for a commutative algebra $R$, and $\Lambda_{R}^{vert, n}(\frg)$, $\Lambda_{R}^{vert, n}(\frg^*)$
are the homogeneous subspaces of degree $n$. 

\subsection{Poisson BRST reduction}
\label{sec:Poisson-brst}

In Sections \ref{sec:Poisson-brst}--\ref{sec:zero-Poisson-BRST}, we construct the jet
bundle of the hypertoric variety $\tldX$ in terms of BRST reduction. The construction
is based on the construction of jet bundles of Slodowy varieties by the BRST reduction
in \cite{AKM}.

Recall that we have the moment map $\mu: T^*V \longrightarrow \frg^*$ and semistable locus
$\frX \subset T^*V$ associated with the torus $G = (\C^\times)^M$-action on the
symplectic vector space $T^*V = T^* \C^N$. Here we took the stability parameter $\delta$
such that the Hamiltonian reduction $X = X_{\delta}$ is a smooth symplectic manifold.
Set $\tldfrX = \frX \times \frg^*$ and $\tldmu: \tldfrX \longrightarrow \frg^*$ as in
\refsec{sec:sympl-deform}. Also, let $\tldX$ be the hypertoric variety as we
introduced in \refsec{sec:sympl-deform}.
The jet bundle
$\calO_{J_{\infty} \tldfrX}$ on $\tldfrX$ is a sheaf of vertex Poisson algebras. By applying the
jet scheme functor $J_{\infty}$ to the moment map $\tldmu : \tldfrX \longrightarrow \frg^*$, we
have a morphism $\tldmu_{\infty} : J_{\infty} \tldfrX \longrightarrow J_{\infty} \frg^*$ and
hence a homomorphism of vertex Poisson algebras
\[
 \tldmu^*_{\infty} : S(\frg \otimes \C[t^{-1}] t^{-1}) \longrightarrow 
 \calO_{J_{\infty} \tldfrX}(\tldfrX)
\]
where the symmetric algebra $S(\frg \otimes \C[t^{-1}] t^{-1})$ has trivial Poisson structure
$Y_{+}(\blkbar, z)$. The homomorphism $\tldmu^*_{\infty}$ is explicitly given by 
$\tldmu^*_{\infty}(A_i) = \sum_{j=1}^{N} \Delta_{ij} x_{j (-1)} y_j - c_i$
for $i=1$, $\dots$, $M$,
where $\Delta = (\Delta_{ij})_{ij}$ is the matrix defined in \refsec{sec:hamilton-act}.

Consider the sheaf of tensor product vertex Poisson algebras
$C_{\VPA} = \calO_{J_{\infty} \tldfrX} \otimes_{\C} \Lambda^{vert}(\frg \oplus \frg^*)$. 
The $\Z$-grading of $\Lambda^{vert}(\frg \oplus \frg^*)$ induces a $\Z$-grading on $C_{\VPA}$
\[
 C_{\VPA} = \bigoplus_{n \in \Z} C_{\VPA}^n, \qquad
 C^n_{\VPA} = \calO_{J_{\infty} \tldfrX} \otimes_{\C} \Lambda^{vert, n}(\frg \oplus \frg^*).
\]
Set $Q_{\VPA} = \sum_{i=1}^{M} \tldmu^*_{\infty}(A_i)_{(-1)} \psi^*_i \in C^1_{\VPA}(\tldfrX)$,
an odd element of degree $+1$ in $C_{\VPA}$. Let 
$d_{\VPA} = Q_{\VPA (0)} = \sum_{i=1}^{M} \sum_{n \in \Z} \tldmu^*_{\infty}(A_i)_{(-n-1)} \psi^*_{i (n)}$
be an operator on $C_{\VPA}$. By definition, the operator $d_{\VPA}$ is a derivation on $C_{\VPA}$. 

\begin{proposition}
 \label{prop:Poisson-coboundary}
 We have $(d_{\VPA})^2 = 0$, and hence, for any open subset $\tldfrU \subset \tldfrX$,
 $(C_{\VPA}(\tldfrU) = \calO_{J_{\infty} \tldfrX} (\tldfrU) \otimes \Lambda^{vert}(\frg \oplus \frg^*), d_{\VPA})$ 
 is a cochain complex.
\begin{proof}
 Since $\tldmu^*_{\infty}$ is a homomorphism of vertex Poisson algebras and the vertex Poisson
 algebra $S(\frg \otimes \C[t^{-1}] t^{-1})$ has trivial Poisson structure, we have 
 $\tldmu^{*}_{\infty}(A_i)_{(n)} \tldmu^*_{\infty}(A_j) = 0$ for any $n \ge 0$ and 
$i$, $j=1$, $\dots$,  $M$. Thus, we have $Q_{\VPA (0)} Q_{\VPA} = 0$. Then, by the axiom of 
vertex Poisson algebras, we have $Q_{\VPA (0)}^2 = (1/2) (Q_{\VPA (0)} Q_{\VPA})_{(0)} = 0$. 
\end{proof}
\end{proposition}

Now we define the notion of the BRST cohomologies for vertex Poisson algebras. 
Take an open subset $\tldfrU \subset \tldfrX$, we consider the cochain complex 
$(C_{\VPA}(\tldfrU), d_{\VPA})$, called a (Poisson) BRST complex. Then, we denote 
its cohomology group
\[
H^{\infty/2 + \bullet}_{\VPA}(\frg, \calO_{J_{\infty} \tldfrX}(\tldfrU)) = H^{\bullet}(C_{\VPA}(\tldfrU), d_{\VPA}), 
\]
and call it (Poisson) BRST cohomology groups.

Note that we have $\partial \circ d_{\VPA} = d_{\VPA} \circ \partial$ following from
$[\partial, Y_+(Q_{\VPA}, z)] = \partial_z Y_+(Q_{\VPA}, z)$. This implies that 
translation operator $\partial$ preserves the subspaces $\Ker d_{\VPA}$ and 
$\Im d_{\VPA} \subset C_{\VPA}(\tldfrU)$. Moreover, by the axiom of vertex Poisson algebras,
the coboundary operator $d_{\VPA} = Q_{\VPA (0)}$ is a derivation with respect to 
${}_{(n)}$-products for all $n \in \Z$. Hence, the $0$-th BRST cohomology 
$H^{\infty/2 + 0}_{\VPA}(\frg, \calO_{J_{\infty} \tldfrX}(\tldfrU)) = H^0(C_{\VPA}(\tldfrU), d_{\VPA})$
is again a vertex Poisson algebra.

Next, we define the BRST cohomologies as a sheaf on the hypertoric variety $\tldX$. 
For an open subset $\tldU \subset \tldX$, let $\tldfrU$ be an open subset of $\tldfrX$ such that 
$\tldfrU$ is closed under the $G$-action and $p^{-1}(\tldU) = \tldfrU \cap \tldmu^{-1}(0)$. The following
lemma asserts that the BRST cohomology 
$H^{\infty/2+\bullet}_{\VPA}(\frg, \calO_{J_{\infty} \tldfrX}(\tldfrU))$ is supported on 
$\tldmu^{-1}(0) \cap \tldfrU$ and it does not depend on the choice of $\tldfrU$. 
Then, we define a sheaf $\calH^{\infty/2 + \bullet}_{\VPA}(\frg, \calO_{J_{\infty} \tldfrX})$
over the hypertoric variety $\tldX$ as the
sheaf associated with the presheaf 
$\tldU \mapsto H^{\infty/2 + \bullet}_{\VPA}(\frg, \calO_{J_{\infty} \tldfrX}(\tldfrU))$ for 
$\bullet \in \Z$.

\begin{lemma}[\cite{AKM}, Thoerem 2.3.2.1]
 \label{lemma:2}
 The presheaf $\tldfrU \mapsto H^{\infty/2 + \bullet}_{\VPA}(\frg, \calO_{J_{\infty} \tldfrX}(\tldfrU))$ 
 over $\tldfrX$ is supported on $\tldmu^{-1}(0)$ and hence it does not depend on the
 choice of $\tldfrU$. 
\end{lemma}

The lemma will be proved in \refsec{sec:double-complex}.

\subsection{Double complex associated with the BRST complex}
\label{sec:double-complex}

The BRST cochain complex can be decomposed into a double cochain complex as follows.
Set
\[
 C_{\VPA}^{p,q} = \calO_{J_{\infty} \tldfrX} \otimes_{\C} \Lambda^{vert, p}(\frg^*) \otimes_{\C} 
 \Lambda^{vert, q}(\frg)
\]
for $p$, $q \in \Z$. Then, we have 
$C_{\VPA}^{n} = \bigoplus_{p + q = n} C_{\VPA}^{p, q}$ for any $n \in \Z$. 
Note that we have 
$\Lambda^{vert, p}(\frg^*) = 0$ unless $p \ge 0$ and $\Lambda^{vert, q}(\frg) = 0$
unless $q \le 0$. Consider the operators
$d^{+}_{\VPA} = \sum_{i=1}^{M} \sum_{n \ge 0} \psi^*_{i (-n-1)} \tldmu^*_{\infty}(A_i)_{(n)}$
and $d^{-}_{\VPA} = \sum_{i=1}^{M} \sum_{n \ge 0} \tldmu^*_{\infty}(A_i)_{(-n-1)} \psi^*_{i (n)}$
on $C_{\VPA}$. Then, $d^{+}_{\VPA}$ maps from $C_{\VPA}^{p, q}$ to $C_{\VPA}^{p+1, q}$,
$d^{-}_{\VPA}$ maps from $C_{\VPA}^{p, q}$ to $C_{\VPA}^{p, q+1}$ and 
we have $d_{\VPA} = d^+_{\VPA} + d^-_{\VPA}$, 
$d^+_{\VPA} \circ d^-_{\VPA} = - d^-_{\VPA} \circ d^+_{\VPA}$. 
Thus, we have a double complex 
$(C_{\VPA}, d^+_{\VPA}, d^-_{\VPA})$ whose total complex is the BRST complex
$(C_{\VPA}, d_{\VPA})$. 

Fix an arbitrary $p \in \Z_{\ge 0}$ and an open subset $\tldfrU \subset \tldfrX$.
Consider the complex 
$(C_{\VPA}^{p, \bullet}(\tldfrU), d^-_{\VPA})$. By the explicit description 
$d^{-}_{\VPA} = \sum_{i=1}^{M} \sum_{n \ge 0} \tldmu^*_{\infty}(A_i)_{(-n-1)} \psi^*_{i (n)}$
of the coboundary operator, the complex $(C_{\VPA}^{p, \bullet}(\tldfrU), d^-_{\VPA})$ coincides
with the Koszul complex of $\calO_{J_{\infty} \tldfrX}(\tldfrU)$ with respect to the sequence
$\{ \tldmu^*_{\infty}(A_i)_{(-n-1)} \}_{i=1, \dots, M, n = 0, 1, \dots}$ (with reversing the degree
of the complex). 
Clearly the sequence $\{ \tldmu^*(A_i) \}_{i=1, \dots, M}$ is a regular sequence 
in $\calO_{\tldfrX}(\tldfrU)$. Then, by the same argument of the proof 
of \cite[Theorem~2.3.3.1]{AKM},
$\{ \tldmu^{*}_{\infty}(A_i)_{(-n-1)} \}_{i, n}$ is also a regular sequence in 
$\calO_{J_{\infty} \tldfrX}(\tldfrU)$. This implies that the cohomology 
$H^{q}(C_{\VPA}^{p, \bullet}(\tldfrU), d_{\VPA}^-)$ vanishes if $q \ne 0$. 
Moreover, when $\tldfrU$ is affine,  we have
$H^{0}(C_{\VPA}^{p, \bullet}(\tldfrU), d_{\VPA}^-) \simeq \calO_{\tldmu^{-1}_{\infty}(0)}(\tldfrU)$
if $\tldmu^{-1}(0) \cap \tldfrU \ne \emptyset$, and zero otherwise for any $p \ge 0$.

Consider the column filtration $\tau_{\ge \bullet} C_{\VPA}(\tldfrU)$; i.e. for $p \in \Z_{\ge 0}$,
$\tau_{\ge p} C_{\VPA}(\tldfrU) = \bigoplus_{k \ge p, q \le 0} C^{k,q}_{\VPA}(\tldfrU)$.
We consider the spectral sequence $\dE_r^{p,q}(\tldfrU)$ associated with the column filtration.
Then we have
$\dE_2^{p,q}(\tldfrU) = H^p(H^q(C_{\VPA}(\tldfrU), d_{\VPA}^-), d_{\VPA}^+)$.

\begin{lemma}
 \label{lemma:conv-spec-Poisson}
 The spectral sequence $\dE_r^{p,q}(\tldfrU)$ converges to the total cohomology
\[
 \dE_r^{p,q}(\tldfrU) \Longrightarrow
 H^{p+q}(C_{\VPA}(\tldfrU), d_{\VPA}) = 
H^{\infty/2 + p+q}_{\VPA}(\frg, \calO_{J_{\infty} \tldfrX}(\tldfrU)).
\]
\begin{proof}
 To prove the convergence, we consider subcomplexes which are bounded both above and below.
 For $m \in \Z_{\ge 0}$, let
 $(C_{\VPA})_m(\tldfrU) = \partial^m (\calO_{\tldfrX}(\tldfrU) \otimes \Lambda(\frg \oplus \frg^*))$
 where we consider
 $\calO_{\tldfrX}(\tldfrU)$ (resp. $\Lambda(\frg \oplus \frg^*)$) as a subalgebra of 
 $\calO_{J_{\infty} \tldfrX}(\tldfrU)$ (resp. $\Lambda^{vert}(\frg \oplus \frg^*)$). Set 
 $(C^{p, q}_{\VPA})_m(\tldfrU) = (C_{\VPA})_m(\tldfrU) \cap C_{\VPA}^{p, q}$. Then, we have
 $C^{p, q}_{\VPA} = \bigoplus_{m \ge 0} (C^{p, q}_{\VPA})_m$.
 By direct computation, for $a \in \calO_{\tldfrX}(\tldfrU)$ and 
 $\varphi \in \Lambda(\frg \oplus \frg^*)$, we have
 \[
  d_{\VPA}(a \otimes \varphi) = \sum_{i=1}^M \{ \tldmu^*(A_i), a \} \otimes \psi^*_{i (-1)} \varphi
 + \sum_{i=1}^{M} \tldmu^*(A_i) a \otimes \{\psi^*_{i (-1)}, \varphi\},
 \]
 and hence $d_{\VPA}$ preserves the subspace $(C_{\VPA})_0(\tldfrU)$. Since 
 $d_{\VPA} = Q_{\VPA (0)}$ commutes with the translation operator $\partial$ by the axiom of
 vertex Poisson algebras, $d_{\VPA}$ also
 preserves $(C_{\VPA})_m(\tldfrU)$ for any $m \in \Z_{\ge 0}$.
 Therefore, 
 $((C_{\VPA})_m(\tldfrU), d^+_{\VPA}, d^-_{\VPA})$ is a double subcomplex of 
 $(C_{\VPA}(\tldfrU), d^+_{\VPA}, d^-_{\VPA})$. 
 Consider the spectral sequence $(E_r^{p,q})_m(\tldfrU)$ associated with the double complex 
 $((C_{\VPA})_m(\tldfrU), d^+_{\VPA}, d^-_{\VPA})$. Since 
$(C_{\VPA})_m(\tldfrU)$  is bounded, the spectral sequence $(E_r^{p,q})_m(\tldfrU)$
converges. This implies the convergence of the spectral sequence $E_r^{p,q}(\tldfrU)$.
\end{proof}
\end{lemma}

As a consequence, we have the BRST cohomology 
$H^{\infty/2+n}_{\VPA}(\frg, \calO_{J_{\infty}\tldfrX}(\tldfrU)) = 0$ for all $n \in \Z$
if $\tldfrU \cap \tldmu^{-1}(0) = \emptyset$.
This implies \reflemma{lemma:2}. Also we have the following vanishing of the negative
BRST cohomologies.

\begin{lemma}
 \label{lemma:vanish-neg-VPA}
 The BRST cohomology
 $H^{\infty/2+n}_{\VPA}(\frg, \calO_{J_{\infty}\tldfrX}(\tldfrU)) = H^{n}(C_{\VPA}(\tldfrU), d_{\VPA})$ 
vanishes if $n < 0$ for any open subset $\tldfrU \subset \tldfrX$. 
\end{lemma}

\subsection{Zeroth Poisson BRST cohomology}
\label{sec:zero-Poisson-BRST}

Now we determine the $0$-th BRST cohomology 
$\calH^{\infty/2+0}_{\VPA}(\frg, \calO_{J_{\infty}\tldfrX})$.
We consider the affine open subset $\tldfrU_J \subset \tldfrX$ introduced in 
\refsec{sec:sympl-deform},
on which the $G$-torsor $\tldmu^{-1}(0) \cap \tldfrU_J \longrightarrow \tldU_J$ is trivial. 
Namely, we have an isomorphism
$\tldmu^{-1}(0) \cap \tldfrU_J \simeq \tldU_J \times G \times \frg^*$ given by the explicit
local coordinate \refeq{eq:tld-triv-G-torsor}.
By applying the functor $J_{\infty}$ to \refeq{eq:tld-triv-G-torsor}, we have 
\begin{multline*}
 \calO_{\tldmu^{-1}_{\infty}(0)}(\tldfrU_J) = 
\C\bigl[\, a^{J*}_{j (-n)}, a^{J}_{j (-n)} \;\bigm\vert\; \substack{j \not\in J \\ n \in \Z_{\ge 1}}\,\bigr] \\
 \otimes_{\C} \C\bigl[\,(T^J_{i (-1)})^{\pm}, T^J_{i (-n-1)} \;\bigm\vert\; \substack{i=1, \dots, M \\ n \in \Z_{\ge 1}}\,\bigr] 
 \otimes_{\C} \C\bigl[\,c_{i (-n)} \;\bigm\vert\; \substack{i=1, \dots, M \\ n \in \Z_{\ge 1}}\,\bigr] \bfone
\end{multline*}
because $\tldmu^{-1}_{\infty}(0) \simeq J_{\infty}(\tldmu^{-1}(0))$ by definition. The action of
$\tldmu^*_{\infty}(A_i)_{(n)} = \gamma_{i (n)}$ for $i=1$, $\dots$, $M$, $n \in \Z_{\ge 0}$,
in the above local coordinate is explicitly given by
$\tldmu^{*}_{\infty}(A_i)_{(n)} = \sum_{k \ge 1} T_{i (-k)} \partial / \partial T_{i (-n-k)}$
by direct calculation. Note that this action coincides with the action of $\frg[t]$ induced 
from the regular representation of $J_{\infty} G$ on 
$\C[J_{\infty} G] \subset \calO_{\tldmu^{-1}(0)}(\tldfrU_J)$.
Since $\dE_1^{p, q}(\tldfrU_J) \simeq \calO_{\tldmu_{\infty}^{-1}(0)}(\tldfrU_J)$
if $q=0$ and zero otherwise, we have $\dE^{0,0}_2(\tldfrU_J) = \Ker d^+_{\VPA}$ where 
\[
d^+_{\VPA} = \sum_{i=1}^{M} \sum_{n \ge 0} \psi^{*}_{i (-n-1)} \tldmu^{*}_{\infty}(A_i)_{(n)} = \sum_{i=1}^{M} \sum_{n,k \ge 0} \psi^*_{i (-n-1)} T_{i (-k)} \partial / \partial T_{i (-n-k)} 
\]
in the above local coordinate. Thus, we have 
\[
 \dE^{0,0}_2(\tldfrU_J) \simeq
\C\bigl[\, a^{*J}_{j (-n)}, a^{J}_{j (-n)} \;\bigm\vert\; \substack{j \not\in J \\ n \in \Z_{\ge 1}}\,\bigr] \\
 \otimes_{\C} \C\bigl[\,c_{i (-n)} \;\bigm\vert\; \substack{i=1, \dots, M \\ n \in \Z_{\ge 1}}\,\bigr] \bfone
\simeq \calO_{J_{\infty} \tldX}(\tldU_J)
\]
and $\dE_r^{0,0}(\tldfrU_J)$ collapses at $r=2$. Therefore, we have
\begin{equation}
\label{eq:1}
 H^{\infty/2+0}_{\VPA}(\frg, \calO_{J_{\infty}\tldfrX}(\tldfrU_J)) \simeq 
 \calO_{J_{\infty} \tldX}(\tldU_J)
\end{equation}
by \reflemma{lemma:conv-spec-Poisson}.

We have the affine open covering $\tldX = \bigcup_J \tldU_J$;
For each indices $I$ and $J$, 
we have the coordinate transformation of 
$\tldvarphi^*_{IJ} : \calO_{\tldfrU_I} \longrightarrow \calO_{\tldfrU_J}$
introduced in \refsec{sec:sympl-deform}.
Its restriction gives the coordinate transformation
$\tldvarphi^*_{IJ} : \calO_{\tldU_I} \longrightarrow \calO_{\tldU_J}$. Applying the
jet scheme functor $J_{\infty}$, we have the isomorphisms
$J_{\infty} \tldvarphi^*_{IJ} : \calO_{J_\infty \tldfrU_I} \vert_{\tldfrU_I \cap \tldfrU_J} \longrightarrow \calO_{J_{\infty} \tldfrU_J} \vert_{\tldfrU_I \cap \tldfrU_J}$
and
$J_{\infty} \tldvarphi^*_{IJ} : \calO_{J_\infty \tldU_I} \vert_{\tldU_I \cap \tldU_J} \longrightarrow \calO_{J_{\infty} \tldU_J} \vert_{\tldU_I \cap \tldU_J}$.
These coordinate transformations are compatible with the isomorphism \refeq{eq:1}, and thus we have 
the following isomorphism of sheaves of Poisson algebras:
\begin{equation}
 \label{eq:Poisson-BRST-isom}
  \calH^{\infty/2+0}_{\VPA}(\frg, \calO_{J_{\infty} \tldfrX}) \simeq \calO_{J_{\infty} \tldX}
\end{equation}
by gluing up $\{\calO_{J_{\infty} \tldU_J}\}_J$ with $\{J_{\infty} \tldvarphi^*_{IJ}\}_{I,J}$.

In the rest of this section, we discuss the BRST reduction 
$H^{\infty/2+0}_{\VPA}(\frg, \calO_{J_{\infty} \tldfrX}(\tldfrX))$ of the coordinate
ring $\calO_{J_{\infty} \tldfrX}(\tldfrX) \simeq \C[J_{\infty} (T^* \C^N \times \frg^*)]$.
Recall the decomposition of the BRST complex, which is introduced in the proof of 
\reflemma{lemma:conv-spec-Poisson},
$(C_{\VPA} = \bigoplus_{m \ge 0} (C_{\VPA})_m, d_{\VPA})$ where 
$(C_{\VPA})_0 = \calO_{\tldfrX} \otimes_{\C} \Lambda(\frg \oplus \frg^*)$ and
$(C_{\VPA})_m = \partial^m (C_{\VPA})_0$ for $m \ge 1$. The subcomplex 
$((C_{\VPA})_0 = \calO_{\tldfrX} \otimes_{\C} \Lambda(\frg \oplus \frg^*), d_{\VPA})$
coincides with the Poisson BRST complex of the Poisson algebra $\calO_{\tldfrX}$ by 
the comoment map 
$\mu^* : \frg \longrightarrow \calO_{\tldfrX}(\tldfrX) \simeq \C[T^*\C^N \times \frg^*]$.
For the detail of the fundamental properties of BRST cohomology of associative algebras,
refer \cite{BRST}. By similar arguments to the above (see also \cite[Section 6.3]{BRST}), we have 
$H^0((C_{\VPA})_0(\tldfrU), d_{\VPA}) \simeq (\calO_{\tldfrX}(\tldfrU) / \sum_{i=1}^M \calO_{\tldfrX}(\tldfrU) \tldmu^*(A_i))^G$
for any open subset $\tldfrU \subset \tldfrX$. Let 
$\calH^0((C_{\VPA})_0, d_{\VPA})$ be the sheaf over $\tldX$ associated with the presheaf
$\tldU \mapsto H^0((C_{\VPA})_0(\tldfrU), d_{\VPA})$ where we take $\tldfrU$ --- an open
subset of $\tldfrX$ which is preserved by the action of $G$ and 
$(\tldmu^{-1}(0) \cap \tldfrU) / G = \tldU$.
Then, we have $\calH^0((C_{\VPA})_0, d_{\VPA}) \simeq \calO_{\tldX}$ and 
$\Gamma(\tldX, \calH^0((C_{\VPA})_0, d_{\VPA})) \simeq \calO_{\tldX}(\tldX) \simeq H^0((C_{\VPA})_0(\tldfrX), d_{\VPA})$
because
$\tldX \longrightarrow \tldX_0 = \tldmu^{-1}(0) \qquot G \simeq \Spec[H^0((C_{\VPA})_0(\tldfrX), d_{\VPA})]$ 
is a resolution of normal singularity.
Since the translation operator $\partial$ commutes with the coboundary operator $d_{\VPA}$,
we have 
$H^0((C_{\VPA})_m(\tldfrX), d_{\VPA}) = \partial^m H^0((C_{\VPA})_0(\tldfrX), d_{\VPA}) \simeq \partial^m \calO_{\tldX}(\tldX)$
for any $m \in \Z_{\ge 0}$. Therefore we have the following proposition.

\begin{proposition}
 \label{prop:global-Poisson-BRST}
 We have 
\[
 H^{\infty/2+0}_{\VPA}(\frg, \calO_{J_{\infty} \tldfrX}(\tldfrX)) = H^0(C_{\VPA}(\tldfrX), d_{\VPA}) \simeq \calO_{J_{\infty} \tldX}(\tldX).
\]
 That is, the Poisson BRST reduction commutes with the global section functor $\Gamma$.
\end{proposition}

\subsection{BRST cohomologies}
\label{sec:def-brst}

Let $\calD^{ch}_{T^* V, \hbar}$ be the sheaf of $\hbar$-adic $\beta\gamma$-system sheaf over the
symplectic vector space $T^* V$ which we defined in \refsec{sec:h-adic-betagamma}. 
By restriction, we define
$\calD^{ch}_{\frX, \hbar} = \calD^{ch}_{T^*V, \hbar} \vert_{\frX}$, the sheaf of 
$\hbar$-adic vertex algebras over $\frX$. Let $V_{\frg, \langle , \rangle, \hbar}$
be a Heisenberg vertex algebra generated by elements $c_1$, $\dots$, $c_M \in \frg$
with the inner product given by
$\langle c_i, c_j \rangle = \sum_{k=1}^{N} \Delta_{ik} \Delta_{jk}$ for $i$, $j=1$, 
$\dots$, $M$. That is, it is the localization over $\frg^*$ of 
the $\hbar$-adic vertex algebra $V_{\frg, \langle, \rangle, \hbar}(\frg^*)$ given by
\[
 V_{\frg, \langle, \rangle, \hbar}(\frg^*) = 
\C[[\hbar]][c_{1 (-n)}, \dots, c_{M (-n)} \;\vert\; n \in \Z_{\ge 1}] \bfone
\]
as a $\C[[\hbar]]$-module, and $c_1$, $\dots$, $c_M$ are bosonic elements whose
 OPEs are given by
$c_i(z) c_j(w) \sim \sum_{k=1}^{N} \hbar^2 \langle c_i, c_j \rangle / (z-w)^2$ where 
$c_i(z) = Y(c_{i (-1)} \bfone, z)$. Set 
$\tldcalD^{ch}_{\frX, \hbar} = \calD^{ch}_{\frX, \hbar} \,\widehat{\otimes}\, V_{\frg, \langle, \rangle, \hbar}$, a sheaf of $\hbar$-adic vertex algebras over $\tldfrX = \frX \times \frg^*$.
Here $\widehat{\otimes}$ is the completion of the tensor product $\otimes_{\C[[\hbar]]}$
with respect to the $\hbar$-adic topology as in \refsec{sec:Clifford}.

To construct the BRST reduction for $\widetilde{\calD}^{ch}_{\frX, \hbar}$, we need to
introduce a quantization of the comoment map $\tldmu^*_{\infty}$. Consider a commutative vertex
algebra $V_0(\frg) = \C[A_{1 (-n)}, \dots, A_{M (-n)}\;\vert\; n \ge 1]$.
Define a $\C[\partial]$-module homomorphism
\begin{equation}
\label{eq:chiral-comoment}
 \mu_{ch}: V_0(\frg) \longrightarrow \widetilde{\calD}^{ch}_{\frX, \hbar}, \quad
  \mu_{ch}(A_i) = \sum_{j=1}^{N} \Delta_{ij} x_{j (-1)} y_j - c_i.
\end{equation}

\begin{lemma}
\label{lemma:chiral-comoment-comm}
 The above map $\mu_{ch}$ preserves the OPEs; i.e. we have
 $\mu_{ch}(A_i)(z) \mu_{ch}(A_j)(w) \sim 0$ for $i$, $j=1$, $\dots$, $M$.
\end{lemma}

This lemma is obviously checked by direct computation.
We call the map $\mu_{ch}$ a chiral comoment map with respect to the $G$-action on
$\frX$.

Consider the sheaf of $\hbar$-adic vertex superalgebras
$C_{\VA} = \widetilde{\calD}^{ch}_{\frX, \hbar} \,\widehat{\otimes}\, Cl^{vert}(\frg \oplus \frg^*)$
where $Cl^{vert}(\frg \oplus \frg^*)$ is the Clifford $\hbar$-adic vertex superalgebra defined in
\refsec{sec:Clifford}. The $\Z$-grading of 
$Cl^{vert}(\frg \oplus \frg^*)$ 
induces a $\Z$-grading on $C_{\VA} = \bigoplus^{\wedge}_{n \in \Z} C_{\VA}^{n}$ where 
$C_{\VA}^{n} = \widetilde{\calD}^{ch}_{\frX, \hbar} \,\widehat{\otimes}\, Cl^{vert,n}(\frg \oplus \frg^*)$
and $\bigoplus^{\wedge}$ is the completion of the direct sum with respect to the $\hbar$-adic
topology.
Consider an odd element $Q_{\VA} = \sum_{i=1}^{M} \mu_{ch}(A_i)_{(-1)} \psi^*_i$ of
degree $+1$ in $C_{\VA}$. Note that the image of $Q_{\VA (0)}$ lies in 
$\hbar \,C_{\VA}$.
Let 
$d_{\VA} = \hbar^{-1} Q_{\VA (0)} = \hbar^{-1} \sum_{i=1}^{M} \sum_{n \in \Z} \mu_{ch}(A_i)_{(-n-1)} \psi^*_{i (n)}$
be a derivation on $C_{\VA}$ homogeneous of degree $+1$.

\begin{proposition}
 \label{prop:chiral-coboundary}
 We have $(d_{\VA})^2 = 0$, and hence, for any open subset $\tldfrU \subset \tldfrX$,
 $(C_{\VA}(\frU) = \widetilde{\calD}^{ch}_{\frX, \hbar}(\tldfrU) \,\widehat{\otimes}\, Cl^{vert}(\frg \oplus \frg^*),\, d_{\VA})$
 is a cochain complex.
 \begin{proof}
  We use the same argument in \refprop{prop:Poisson-coboundary}. By 
  \reflemma{lemma:chiral-comoment-comm} we have $\mu_{ch}(A_i)_{(n)} \mu_{ch}(A_j) = 0$
  for all $i$, $j = 1$, $\dots$, $M$ and $n \ge 0$. Thus, we have $Q_{\VA (0)} Q_{\VA} = 0$. 
  By Borcherds' identity, we have 
\[
(Q_{\VA (0)})^2 = (1/2)[Q_{\VA (0)}, Q_{\VA (0)}] = (1/2)(Q_{\VA (0)} Q_{\VA})_{(0)} = 0. 
\]
 \end{proof}
\end{proposition}

Now we define the notion of the chiral BRST cohomologies.
Taking an open subset $\tldfrU \subset \tldfrX$, we consider the cochain complex 
$(C_{\VA}(\tldfrU), d_{\VA})$, called a BRST complex. Then, for $n \in \Z$, we denote 
its cohomology group
$H^{\infty/2 + n}_{\VA}(\frg, \widetilde{\calD}^{ch}_{\frX, \hbar}(\tldfrU)) = H^{n}(C_{\VA}(\tldfrU), d_{\VA})$, and call it the $n$-th BRST cohomology.
Note that we have $[\partial, Y(Q_{\VA}, z)] = \partial_z Y(Q_{\VA}, z)$ on $C_{\VA}$ by the axiom 
of $\hbar$-adic vertex superalgebras. By taking the coefficient of $z^{-1}$, we obtain
$\partial \circ Q_{\VA(0)} - Q_{\VA(0)} \circ \partial = 0$. Thus, the translation operator
$\partial$ preserves the subspaces $\Ker d_{\VA}$ and $\Im d_{\VA}$. Further, for any element
$a$, $b \in C_{\VA}$ and for any $n \in \Z$, we have 
$Q_{\VA(0)} a_{(n)} b - (-1)^{\prty{a}} a_{(n)} Q_{\VA(0)} b = (Q_{\VA(0)} a)_{(n)} b$ by
the Borcherds' identity. By taking $a$, $b$ from $\Ker d_{\VA}$, we conclude that 
$Q_{\VA(0)} (a_{(n)} b) = 0$ and thus $a_{(n)} b \in \Ker d_{\VA}$. Also, by taking 
$a \in \Ker d_{\VA}$ and $b \in C_{\VA}$, we have
$a_{(n)} (Q_{\VA(0)} b) = Q_{\VA(0)}(a_{(n)} b) \in \Im d_{\VA}$. Therefore we conclude the
following proposition.

\begin{proposition}
 \label{prop:BRST-is-VA}
 For an open subset $\tldfrU \subset \tldfrX$,
 the $0$-th BRST cohomology group
 $H^{\infty/2+0}_{\VA}(\frg, \widetilde{\calD}^{ch}_{\frX, \hbar}(\tldfrU)) = H^0(C_{\VA}(\tldfrU), d_{\VA})$
 is an $\hbar$-adic vertex algebra.
\end{proposition}

Next, we define the BRST cohomology group as a sheaf on the hypertoric variety $\tldX$. 
For an open subset $\tldU \subset \tldX$, let $\tldfrU$ be an open subset of $\tldfrX$ such that 
$\tldfrU$ is closed under the $G$-action and 
$p^{-1}(\tldU) = \tldfrU \cap \tldmu^{-1}(0)$. The following
lemma asserts that the BRST cohomology group
$H^{\infty/2+\bullet}_{\VA}(\frg, \widetilde{\calD}^{ch}_{\frX, \hbar}(\tldfrU))$ is supported on 
$\tldmu^{-1}(0) \cap \tldfrU$ and it does not depend on the choice of $\tldfrU$. 
Then, we define a sheaf 
$\calH^{\infty/2 + \bullet}_{\VA}(\frg, \widetilde{\calD}^{ch}_{\frX, \hbar})$
over the hypertoric variety $\tldX$ as the
sheaf associated with the presheaf 
$\tldU \mapsto H^{\infty/2 + \bullet}_{\VA}(\frg, \widetilde{\calD}^{ch}_{\frX, \hbar}(\tldfrU))$.

\begin{lemma}[cf. \cite{AKM}, Thoerem 2.3.5.1]
\label{lemma:1}
 The presheaf 
 $\tldfrU \mapsto H^{\infty/2 + \bullet}_{\VA}(\frg, \widetilde{\calD}^{ch}_{\frX, \hbar}(\tldfrU))$ 
 over $\tldfrX$ is supported on $\tldmu^{-1}(0)$ and hence it does not depend on the
 choice of $\tldfrU$. 
\end{lemma}

In the rest of this section, we prove \reflemma{lemma:1}.
The coboundary operator of the BRST complex 
$d_{\VA} = \hbar^{-1} Q_{\VA (0)} = \hbar^{-1} \sum_{i=1}^{M} \sum_{n \in \Z} \mu_{ch}(A_i)_{(-n-1)} \psi^*_{i(n)}$
is separated into two parts $d^+_{\VA}$ and $d^-_{\VA}$; namely, putting 
\[
d^+_{\VA} = \hbar^{-1} \sum_{i=1}^{M} \sum_{n \ge 0} \psi^*_{i (-n-1)} \mu_{ch}(A_i)_{(n)},
\quad
d^-_{\VA} = \hbar^{-1} \sum_{i=1}^{M} \sum_{n \ge 0} \mu_{ch}(A_i)_{(-n-1)} \psi^*_{i (n)}, 
\]
we have $d_{\VA} = d^+_{\VA} + d^-_{\VA}$. 
Moreover, we have 
$d^+_{\VA} \circ d^-_{\VA} = - d^-_{\VA} \circ d^+_{\VA}$ because 
$\mu_{ch}(A_i)_{(n)} \mu_{ch}(A_j) = 0$ and 
$\psi^*_{i (n)} \psi^*_{j} = 0$ for any $i$, $j = 1$, $\dots$, $M$ and $n \ge 0$.
Thus, we have a double complex $(C_{\VA}, d^+_{\VA}, d^-_{\VA})$
where 
\[
 C_{\VA}^{p,q} = \tldcalD^{ch}_{\frX, \hbar} \,\widehat{\otimes}\, \Lambda^{vert, p}_{\C[[\hbar]]}(\frg^*) \,\widehat{\otimes}\, \Lambda^{vert, q}_{\C[[\hbar]]}(\frg)
\]
for $p$, $q \in \Z$, is induced from the decomposition of the $\hbar$-adic Clifford vertex 
algebra. Note that $C_{\VA}^{p,q} = 0$ unless $p \ge 0$ and $q \le 0$; that is,
$(C_{\VA}, d^+_{\VA}, d^-_{\VA})$ is the fourth quadrant cochain double 
complex.

The BRST complex 
$C_{\VA} = \widetilde{\calD}^{ch}_{\frX, \hbar} \,\widehat{\otimes}\, Cl^{vert}(\frg \oplus \frg^*)$
is naturally equipped with a filtration $F_{\bullet} C_{\VA}$ by powers of $\hbar$: 
$F_p C_{\VA} = \hbar^p C_{\VA}$ for $p \in \Z_{\ge 0}$. 
For each $p \in \Z_{\ge 0}$, the associated graded space is 
$\Gr_{p} C_{\VA} = F_{p} C_{\VA} / F_{p+1} C_{\VA} \simeq \calO_{J_{\infty} \tldfrX} \otimes \Lambda^{vert}(\frg \oplus \frg^*)$
as vertex Poisson superalgebras.

Consider the action of 
$d^+_{\VA}$ and $d^-_{\VA}$ on the vertex Poisson superalgebra
$\Gr_{p} C_{\VA} \simeq \calO_{J_{\infty} \tldfrX} \otimes_{\C} \Lambda^{vert}(\frg \oplus \frg^*)$. 
The operators $d^+_{\VA}$, $d^-_{\VA}$ act by 
\[
\bar{d}^+_{\VA} = \sum_{i=1}^{M} \sum_{n \ge 0} \psi^*_{i (-n-1)} \tldmu_{\infty}^*(A_i)_{(n)}, 
\quad
\bar{d}^-_{\VA} = \sum_{i=1}^{M} \sum_{n \ge 0} \tldmu^*_{\infty}(A_i)_{(-n-1)}  \psi^*_{i (n)}
\]
respectively on $\calO_{J_{\infty}\tldfrX} \otimes_{\C} \Lambda^{vert}(\frg \oplus \frg^*)$.

Thus, for each $p \in \Z_{\ge 0}$, the double complex 
$(\Gr_p C_{\VA}, \bar{d}^+_{\VA}, \bar{d}^-_{\VA})$ is isomorphic to
the double complex $(C_{\VPA}, d^+_{\VPA}, d^-_{\VPA})$ associated with the
Poisson BRST complex which we discussed in \refsec{sec:double-complex}.
By \reflemma{lemma:2}, for an open subset $\tldfrU \subset \tldfrX$ such that 
$\tldfrU \cap \tldmu^{-1}(0) = \emptyset$, we have 
\[
 H^{\bullet}(\Gr_p C_{\VA}(\tldfrU), \bar{d}^+_{\VA} + \bar{d}^-_{\VA}) = 
 H^{\bullet}(C_{\VPA}(\tldfrU), d_{\VPA}) = 0
\]
for any $\bullet \in \Z$ and $p \in \Z_{\ge 0}$. Now we consider the spectral sequence
$E^{p,q}_r$ associated with the filtered complex $(F_{\bullet} C_{\VA}, d_{\VA})$. 
Then, we have $E_1^{p,q}(\tldfrU) = 0$ for any $p$, $q \in \Z$ by the above and thus 
$E_r^{p,q}(\tldfrU)$
collapses at $r=1$. Since the filtration $F_{\bullet} C_{\VA}$ is bounded above and complete,
$E_r^{p,q}(\tldfrU)$ converges to $\Gr_p H^{p+q}(C_{\VA}(\tldfrU), d_{\VA})$ by the complete 
convergence theorem \cite[Theorem 5.5.10]{Weibel}. Therefore, we have 
the vanishing 
$H^{\bullet}(C_{\VA}(\tldfrU), d_{\VA}) = 0$ for an open subset $\tldfrU$ which intersects 
trivially with $\tldmu^{-1}(0)$, and it proves \reflemma{lemma:1}.

By a similar argument, we obtain the vanishing of negative BRST cohomologies as follows.

\begin{proposition}
 \label{prop:vanish-neg-BRST}
 For $n < 0$, we have
 $\calH^{\infty/2 + n}_{\VA}(\frg, \tldcalD^{ch}_{\frX, \hbar}) = 0$.
\begin{proof}
For any $p \ge 0$, $n < 0$ and any open subset $\tldfrU \subset \tldfrX$, we have
\[
 H^{n}(\Gr_p C_{\VA}(\tldfrU), \bar{d}_{\VA}^+ + \bar{d}_{\VA}^-)
= H^{n}(C_{\VPA}(\tldfrU), d_{\VPA}) = 0
\]
by \reflemma{lemma:vanish-neg-VPA}. Again we consider the spectral sequence $E^{p,q}_r$
associated with the $\hbar$-adic filtration. Then, for any $p$, $q$ such that $p+q < 0$,
we have $E^{p,q}_1(\tldfrU) = 0$. Since the filtration is complete and bounded above, 
and we have $E^{p,q}_1(\tldfrU) = E_2^{p,q}(\tldfrU) = \dots$ for $p$, $q$ with
$p+q < 0$, the spectral sequence $E_r^{p,q}(\tldfrU)$ converges to the cohomology
$\Gr_p H^{p+q}(C_{\VA}(\tldfrU), d_{\VA})$ when $p+q < 0$. Thus the cohomology
$H^n(C_{\VA}(\tldfrU), d_{\VA})$ vanishes for negative $n$. 
\end{proof}
\end{proposition}

\begin{definition}
 \label{def:BRST-red-sheaf}
 We write
 the $0$-th cohomology $\calH^{\infty/2+0}_{\VA}(\frg, \tldcalD^{ch}_{\frX, \hbar})$
 by $\tldcalD^{ch}_{X, \hbar}$.
\end{definition}

\subsection{Spectral sequence associated with the double complex}
\label{sec:conv-double-cpx}

For any fixed $k \in \Z$, we have the complex $(C^{k, \bullet}_{\VA}, d^{-}_{\VA})$.
Set a $\C[[\hbar]]$-submodule
\[
 C^{\bullet}_{-} = \tldcalD^{ch}_{\frX, \hbar} \,\widehat{\otimes}\,
\Lambda^{vert, \bullet}_{\C[[\hbar]]}(\frg)
\subset 
\tldcalD^{ch}_{\frX, \hbar} \,\widehat{\otimes}\,
\Lambda^{vert, \bullet}_{\C[[\hbar]]}(\frg) \,\widehat{\otimes}\,
\Lambda^{vert,k}_{\C[[\hbar]]}(\frg^*) = C^{k, \bullet}_{\VA},
\]
and we have a complex $(C_{-}, d^{-}_{\VA})$.
Note that 
\[
 H^{n}(C_{\VA}^{k,\bullet}, d^{-}_{\VA}) \simeq H^{n}(C_{-}^{\bullet}, d^{-}_{\VA}) \,\widehat{\otimes}\, \Lambda^{vert, k}_{\C[[\hbar]]}(\frg^*)
\]
for any $k$, $n \in \Z$.

Consider the filtration of $C_{-}$ given by the powers of $\hbar$,
denoted $F_{p} C_{-} = \hbar^p C_{-}$ ($p \ge 0$).
Clearly, the coboundary operator $d^{-}_{\VA}$ preserves the filtration.
Let $E_{r}^{p,q}$ be the spectral sequence associated with the filtration 
$F_{\bullet} C_{-}$. Then, we have 
$E^{p,q}_0 = F_p C^{p+q}_{-} / F_{p+1} C^{p+q}_{-} \simeq \calO_{J_{\infty} \tldfrX} \otimes_{\C} \Lambda^{vert, p+q}_{\C}(\frg)$ on which
the coboundary operator acts by $\bar{d}^{-}_{\VA} = d^{-}_{\VPA}$ as we see in the
previous section. Take an open subset $\tldfrU \subset \tldfrX$.
By the result of \refsec{sec:double-complex}, we have
\[
 E_1^{p,q}(\tldfrU) \simeq 
\begin{cases}
  \calO'_{\tldmu^{-1}_{\infty}(0)}(\tldfrU) &
 (p+q=0) \\
 0 & (p+q \ne 0)
\end{cases}
\]
where
\[
 \calO'_{\tldmu^{-1}_{\infty}(0)}(\tldfrU) = 
\calO_{J_{\infty}\tldfrX}(\tldfrU) \bigm/ \sum_{i=1}^{M} \sum_{n \ge 0} \tldmu^{*}_{\infty}(A_i)_{(-n-1)} \calO_{J_{\infty} \tldfrX}(\tldfrU).
\]
Note that $\calO'_{\tldmu^{-1}_{\infty}(0)}(\tldfrU) = \calO_{\tldmu^{-1}_{\infty}(0)}(\tldfrU)$
when $\tldfrU$ is an affine open subset.
The above implies that the spectral sequence
$E_{r}^{p,q}(\tldfrU)$ collapses at $r=1$. Since the filtration 
$F_{\bullet} C_{-}(\tldfrU)$ is bounded above and complete, the collapse
implies the convergence of the spectral sequence by the complete convergence theorem 
\cite[Theorem 5.5.10]{Weibel}.

\begin{lemma}
 \label{lemma:BRST-minus-cohom}
 For $p \ge 0$ and an open subset $\tldfrU \subset \tldfrX$, we have an 
 isomorphism
\[
 H^q(C^{p, \bullet}_{\VA}(\tldfrU), d^-_{\VA}) 
 \simeq 
 \calO'_{\tldmu^{-1}_{\infty}(0)}(\tldfrU)[[\hbar]] \,\widehat{\otimes}\, \Lambda^{vert, p}_{\C[[\hbar]]}(\frg^*)
\]
if $q=0$, and zero otherwise.
\end{lemma}

Consider the double complex
$(C_{\VA}(\tldfrU), d^+_{\VA}, d^-_{\VA})$. 
Consider the column filtration $\tau_{\ge \bullet} C_{\VA}$ of the double complex $C_{\VA}$;
$\tau_{\ge p} C_{\VA}(\tldfrU) = \bigoplus^{\wedge}_{k \ge p, q \le 0} C^{k,q}_{\VA}(\tldfrU)$
where $\bigoplus^{\wedge}$ is the completion with respect to the $\hbar$-adic topology.
Let $\dE^{p,q}_{r}(\tldfrU)$ be the spectral sequence
associated with the column filtration $\tau_{\ge \bullet} C_{\VA}$.
By \reflemma{lemma:BRST-minus-cohom}, we have 
\begin{multline*}
 \dE^{p,q}_2(\tldfrU) = H^p(H^q(C_{\VA}(\tldfrU), d^-_{\VA}), d^+_{\VA}) \\
 \simeq
 \begin{cases}
  H^p(\calO'_{\tldmu^{-1}_{\infty}(0)}(\tldfrU)[[\hbar]] \,\widehat{\otimes}\, 
  \Lambda^{vert, \bullet}_{\C[[\hbar]]}(\frg^*), d^{+}_{\VA}) & \text{if } q=0, \\
  0 & \text{otherwise.}
 \end{cases}
\end{multline*}
Thus, the spectral sequence $\dE^{p,q}_r(\tldfrU)$ collapses at $r=2$.
In the rest of this section we prove the following proposition.

\begin{proposition}
 \label{prop:conv-double-cpx}
 The spectral sequence $\dE^{p,q}_r(\tldfrU)$ converges to the total cohomology
 $H^{p+q}(C_{\VA}(\tldfrU), d_{\VA})$.
\end{proposition}

Consider the completion $\widehat{C}_{\VA}(\tldfrU)$ of the BRST complex 
$C_{\VA}(\tldfrU)$ with respect to the column filtration $\tau_{\ge \bullet} C_{\VA}(\tldfrU)$.
Then, the column filtration $\tau_{\ge \bullet} \widehat{C}_{\VA}(\tldfrU)$ is bounded
above, complete, and the spectral sequence $\dE^{p,q}_r(\tldfrU)$ collapses at $r=2$. Thus, 
the spectral sequence $\dE^{p,q}_r(\tldfrU)$ converges to the total cohomology group
$H^{p+q}(\widehat{C}_{\VA}(\tldfrU), d_{\VA})$ of the completed complex
by the complete convergence theorem 
\cite[Theorem 5.5.10]{Weibel}. 

For $p$, $q \in \Z$, set 
\begin{align*}
 A^{p,q}_{\infty} &= 
 \{\, \alpha \in \tau_{\ge p} C_{\VA}^{p+q}(\tldfrU) \;|\; d_{\VA}\alpha = 0 \,\}, \\
 \widehat{A}^{p,q}_{\infty} &=
 \{\, \alpha \in \tau_{\ge p} \widehat{C}_{\VA}^{p+q}(\tldfrU) \;|\; d_{\VA}\alpha = 0 \,\}.
\end{align*}
To prove the convergence of the spectral sequence $\dE^{p,q}_{r}(\tldfrU)$ to 
$H^{p+q}(C_{\VA}(\tldfrU), d_{\VA})$, it is sufficient to show 
$A^{p,q}_{\infty} = \widehat{A}^{p,q}_{\infty}$ and 
$\bigcap_{p \ge 0} \tau_{\ge p} H^{n}(C_{\VA}(\tldfrU), d_{\VA}) = 0$ for any $n$.

\begin{lemma}
 \label{lemma:5}
 We have $A^{p,q}_{\infty} = \widehat{A}^{p,q}_{\infty}$ for any $p$, $q \in \Z$.
\begin{proof}
 Take arbitrary $\alpha = \sum_{k = 0}^{\infty} \alpha_k \in \widehat{A}^{p,q}_{\infty}$ where
 $\alpha_k \in C^{p+k, q-k}_{\VA}$. We show that $\alpha$ belongs to $A^{p,q}_{\infty}$.
 By the condition, we have $d_{\VA}^{-} \alpha_0 = 0$ and 
 $d_{\VA}^{+} \alpha_{k} = - d_{\VA}^{-} \alpha_{k+1}$ for $k \ge 0$.
 Since the coboundary operators $d_{\VA}^{+}$, $d_{\VA}^{-}$ of the double complex 
 preserve the $\hbar$-adic filtration $F_{\bullet} C_{\VA}(\tldfrU)$, we may assume that
 $\alpha_k \in F_{s} C_{\VA}$ implies $\alpha_{k+1} \in F_{s} C_{\VA}(\tldfrU)$ for any 
 $k$, $s \ge 0$.
 Considering modulo $F_1 C_{\VA}(\tldfrU)$, we have the isomorphism of double complexes
 $(C_{\VA} / F_1 C_{\VA}, \bar{d}_{\VA}^{+}, \bar{d}_{\VA}^{-}) \simeq (C_{\VPA}, d_{\VPA}^+, d_{\VPA}^{-})$, 
 which we showed in \refsec{sec:def-brst}. The double complex $C_{\VPA}(\tldfrU)$ is decomposed into
 a direct sum of bounded double complexes, and thus the associated spectral sequence converges
 as \reflemma{lemma:conv-spec-Poisson}. This implies that there exists
 an integer $k_1 \ge 0$ such that $\alpha_l \in F_1 C_{\VA}(\tldfrU)$ for all $l \ge k_1$. 

 Assume that, for an integer $s \ge 0$, there exists 
 $k_s \ge 0$ such that $\alpha_l \in F_s C_{\VA}(\tldfrU)$ for $l \ge k_s$. By the condition,
 we have $d_{\VPA}^+ \overline{\alpha}_k = - d_{\VPA}^{-} \overline{\alpha}_{k+1}$ for 
 $k \ge k_s$, where $\overline{\alpha}_k$ is the image of $\alpha_k$ in 
 $F_s C_{\VA}(\tldfrU) / F_{s+1} C_{\VA}(\tldfrU) \simeq C_{\VPA}(\tldfrU)$. By the above 
 equalities, $\overline{\alpha}_k$
 for $k \ge k_s$ belong to the same bounded double subcomplex, and hence we have
 $\overline{\alpha}_{k_{s+1}} = \overline{\alpha}_{k_{s+1}+1} = \dots = 0$ for some $k_{s+1} \ge 0$.
 Thus, inductively on $s \ge 0$, we have an integer $k_s \ge 0$ such that
 $\alpha_l \in F_s C_{\VA}(\tldfrU)$ for any $l \ge k_s$. This implies that 
 $\alpha \in C_{\VA}(\tldfrU)$,
 and thus $\alpha \in A^{p,q}_{\infty}$.
\end{proof}
\end{lemma}

The above lemma asserts that the spectral sequence $\dE^{p,q}_{r}(\tldfrU)$ weakly converges to
the total cohomology $H^{p+q}(C_{\VA}(\tldfrU), d_{\VA})$; that is we have
\[
\dE^{p,q}_{2}(\tldfrU) = \dots = \dE^{p,q}_{\infty}(\tldfrU) \simeq \Gr_{p} H^{p+q}(C_{\VA}(\tldfrU), d_{\VA}) 
\]
 for any $p$, $q \in \Z$.

\begin{lemma}
 \label{lemma:6}
 For any $n \in \Z$, we have $\bigcap_{p \ge 0} \tau_{\ge p} H^n(C_{\VA}(\tldfrU), d_{\VA}) = 0$.
 \begin{proof}
  Take arbitrary $\alpha \in \bigcap_{p \ge 0} \tau_{\ge p} H^n(C_{\VA}(\tldfrU), d_{\VA})$.
  Let $\sum_{k=0}^{\infty} \alpha_k$ where
  $\alpha_k \in C^{n+k,-k}_{\VA}(\tldfrU)$ be a cocycle which represent $\alpha$. By the
  convergence of $\dE^{p,q}_{r}(\tldfrU)$, we have 
  $\bigcap_{p \ge 0} \tau_{\ge p} H^n(\widehat{C}_{\VA}(\tldfrU), d_{\VA}) = 0$,
  and thus we have $\alpha' = \sum_{k=0}^{\infty} \alpha'_k \in \widehat{C}_{\VA}(\tldfrU)$
  with $\alpha'_k \in C_{\VA}^{n+k, -k}$ such that $d_{\VA} \alpha' = \sum_{k=0}^{\infty} \alpha_k$.
  We have $d_{\VA}^{-} \alpha'_0 = \alpha_0$ and 
  $d_{\VA}^{+} \alpha'_k + d_{\VA}^{-} \alpha'_{k+1} = \alpha_{k+1}$ for $k \ge 0$. 
  Let $\overline{\alpha}_k$ (resp. $\overline{\alpha}'_k$) be the image of $\alpha_k$
  (resp. $\alpha'_k$) in $C_{\VA}(\tldfrU) / F_1 C_{\VA}(\tldfrU) \simeq C_{\VPA}(\tldfrU)$.
  Then, we have equalities 
  $d_{\VPA}^{-} \overline{\alpha}'_0 = \overline{\alpha}_0$ and 
  $d_{\VPA}^{+} \overline{\alpha}'_k + d_{\VPA}^{-} \overline{\alpha}'_{k+1} = \overline{\alpha}_{k+1}$
  for $k \ge 0$. Note that we have $\overline{\alpha}'_k \ne 0$ for finitely many $k$ because
  $\alpha \in C_{\VA}(\tldfrU)$.
  Consider a bounded double subcomplex which contains 
  $\overline{\alpha}'_0$, $\overline{\alpha}_0$, $\dots$. By the above
  equalities, $\overline{\alpha}'_k$ for $k \ge 0$ also belongs to the same bounded 
  double subcomplex. Thus, there exists $k_1 \ge 0$ such that 
  $\overline{\alpha}'_l = 0$ i.e. $\alpha'_l \in F_1 C_{\VA}(\tldfrU)$ for $l \ge k_1$.
  By the same argument of the proof of \reflemma{lemma:5}, inductively on $s$, there exists
  $k_s \ge 0$ such that $\alpha'_l \in F_s C_{\VA}(\tldfrU)$ for $l \ge k_s$.
  Therefore, $\alpha' = \sum_{k=0}^{\infty} \alpha'_k \in C_{\VA}(\tldfrU)$, and hence
  $\alpha = d_{\VA} \alpha' = 0$ in $H^n(C_{\VA}(\tldfrU), d_{\VA})$.
 \end{proof}
\end{lemma}

\reflemma{lemma:6} together with \reflemma{lemma:5} gives the convergence of the spectral
sequence $\dE^{p,q}_{r}(\tldfrU)$ to the total cohomology $H^{p+q}(C_{\VA}(\tldfrU), d_{\VA})$
(\refprop{prop:conv-double-cpx}).

\section{Local structure of BRST reduction}
\label{sec:local-struct}

In the previous sections, we defined the sheaf of $\hbar$-adic vertex algebras 
$\tldcalD^{ch}_{X, \hbar}$ over the hypertoric variety $\tldX$. 
Now we describe the local structure of the BRST reduction
$\tldcalD^{ch}_{X, \hbar}(\tldU_J)$ over 
the affine open subset $\tldU_J \subset \tldX$ with using the local coordinate 
which we defined in \refsec{sec:triv-Hamil-red}.

Consider the affine open subset $\tldfrU_J \subset \tldfrX$ defined in \refsec{sec:sympl-deform},
and recall the local coordinate functions in $\calO_{\tldfrX}(\tldfrU_J)$ of \refeq{eq:tld-local-trivial}.
We identify these coordinate functions with their image in $\calO_{J_{\infty} \tldfrX}(\tldfrU_J)$
and their lifts onto $\tldcalD^{ch}_{\frX, \hbar}(\tldfrU_J)$. 
Then, consider the cochains $a_j^{* J} = a_{j (-1)}^{* J} \bfone$, 
$a_j^{J} = a^{J}_{j (-1)} \bfone \in C^0_{\VA}(\tldfrU_J)$ for $j \not\in J$. 
Since $a^{* J}_j$ and $a_j^{J}$ have none of the factors $x_{k (-1)} y_k$ for $k=1$, $\dots$, $N$
and they are $G$-invariant, we have 
$\mu_{ch}(A_i)_{(n)} a^{* J}_j = \mu_{ch}(A_i)_{(n)} a_j^{J} = 0$ for any $i=1$, $\dots$, $M$,
$j \not\in J$ and $n \in \Z_{\ge 0}$. It implies that 
$d_{\VA} a^{* J}_j = d_{\VA} a_j^{J} = 0$ and thus $a^{* J}_j$, $a_j^{J}$ define elements in 
$H^0(C_{\VA}(\tldfrU_J), d_{\VA}) = \tldcalD^{ch}_{X, \hbar}(\tldU_J)$. We denote these
elements the same notation $a^{* J}_j$ and $a^{J}_j \in \tldcalD^{ch}_{X, \hbar}(\tldU_J)$.
For $i=1$, $\dots$, $M$, we have a cochain $c_i = c_{i (-1)} \bfone \in C^0_{\VA}(\tldfrU_J)$.
By direct calculation, we have 
$d_{\VA} c_i = - \hbar \sum_{j=1}^{M} \langle c_i, c_j \rangle \psi^{*}_{j (-2)} \bfone$,
which is not necessarily zero. Note that, on $\tldfrU_J$, we have a cochain 
$\partial \log T_j^{J} = T^{J}_{j (-2)} (T^{J}_{j (-1)})^{-1} \bfone \in C^{0}_{\VA}(\tldfrU_J)$ 
for $j = 1$, $\dots$, $M$. Again since $T^J_j$ has none of the factors $x_{k (-1)} y_k$ for
$k=1$, $\dots$, $N$, and $T^J_j$ is of weight $\bfe_j$ with respect to the $G$-action, 
we have $\mu_{ch}(A_i)_{(n)} \partial \log T^J_j = \delta_{n1} \delta_{ij} \hbar \bfone$
for $i$, $j=1$, $\dots$, $M$ and $n \in \Z_{\ge 0}$. For $i=1$, $\dots$, $M$, set
a cochain locally defined on $\tldfrU_J$,
\begin{equation}
 b^J_i = c_i + \hbar \sum_{j=1}^{M} \langle c_i, c_j \rangle \partial \log T^J_j
\in C^0_{\VA}(\tldfrU_J).
\end{equation}
Then, we have $d_{\VA} b^J_i = 0$, and thus $b^J_i$ defines an element of 
$b^J_i \in H^0(C_{\VA}(\tldfrU_J), d_{\VA}) = \tldcalD^{ch}_{X, \hbar}(\tldU_J)$ for 
$i=1$, $\dots$, $M$. By \refprop{prop:BRST-is-VA}, 
$H^0(C_{\VA}(\tldfrU_J), d_{\VA}) = \tldcalD^{ch}_{X, \hbar}(\tldU_J)$ is an
$\hbar$-adic vertex algebra. Thus, we have
\begin{equation}
\label{eq:3} 
 \C[[\hbar]][\,a^{* J}_{j (-n)}, a^{J}_{j (-n)}, b^{J}_{i (-n)} \;|\; j \not\in J, i=1, \dots, M, 
 n \in \Z_{\ge 1}\,]
 \subset H^0(C_{\VA}(\tldfrU_J), d_{\VA}).
\end{equation}

By \refprop{prop:conv-double-cpx} and \reflemma{lemma:BRST-minus-cohom}, we have
$H^0(C_{\VA}(\tldfrU_J), d_{\VA}) = \dE_2^{p,q}(\tldfrU_J) = H^0(C_{+}(\tldfrU_J), d_{\VA}^+)$
where 
$C_{+}^{\bullet} = \calO_{\tldmu^{-1}_{\infty}(0)}[[\hbar]] \,\widehat{\otimes}\, \Lambda^{vert, \bullet}_{\C[[\hbar]]}(\frg^*)$.
Since $d^{+}_{\VA} = d_{\VPA}^{+}$ on $\Gr C_{+}$, we have an embedding
\begin{multline*}
 H^0(C_{+}(\tldfrU_J), d_{\VA}^{+}) = \Ker d_{\VA}^{+} \cap C^0_{+}(\tldfrU_J) \\
\subset H^0(\calO_{\tldmu^{-1}_{\infty}(0)}(\tldfrU_J)[[\hbar]] \otimes 
\Lambda^{vert, \bullet}_{\C}(\frg^*), d^+_{\VPA}) \simeq
 \calO_{J_{\infty} \tldX}(\tldU_J)[[\hbar]].
\end{multline*}
By \refeq{eq:1}, we have 
\[
\calO_{J_{\infty} \tldX}(\tldU_J)[[\hbar]] = \C[[\hbar]][\, a^{* J}_{j (-n)}, a^{J}_{j (-n)}, c_{i (-n)}\;|\; j \not\in J, i=1, \dots, M, n \in \Z_{\ge 1}\,],
\]
and thus the $\hbar$-adic vertex subalgebra of \refeq{eq:3} coincides with 
$H^0(C_{\VA}(\tldfrU_J), d_{\VA})$. Here note that the elements $a^{* J}_{j (-n)}$, $a^{J}_{j (-n)}$ and
$b^{J}_{i (-n)}$ for $j \not\in J$, $i=1$, $\dots$, $M$ and $n \in \Z_{\ge 1}$ are algebraically
independent because their images $a^{* J}_{j (-n)}$, $a^{J}_{j (-n)}$, $c_{i (-n)}$ in 
$\calO_{J_{\infty}\tldX}(\tldU_J)[[\hbar]] / \hbar\, \calO_{J_{\infty} \tldX}(\tldU_J)[[\hbar]]$
are algebraically independent.

\begin{proposition}
 \label{prop:local-descr}
 For the affine open subset $\tldU_J \subset \tldX$ defined in \refsec{sec:sympl-deform},
 we have 
\begin{align*}
  \tldcalD^{ch}_{X, \hbar}(\tldU_J) &= H^0(C_{\VA}(\tldfrU_J), d_{\VA}) \\
 &=
 \C[[\hbar]][\,a^{* J}_{j (-n)}, a^{J}_{j (-n)}, b^{J}_{i (-n)} \;|\; j \not\in J, i=1, \dots, M, 
 n \in \Z_{\ge 1}\,] \\
 &\simeq \calD^{ch}(T^* \C^{N-M})_{\hbar} \,\widehat{\otimes}\, V_{\langle, \rangle, \hbar}(\frg).
\end{align*}
\begin{proof}
 The isomorphism as $\C[[\hbar]]$-modules follows from the above discussion. 
 We consider the structure as an $\hbar$-adic vertex algebra. Note that, by the explicit
 construction in \refsec{sec:triv-Hamil-red}, $a^{* J}_{j}$, $a^{J}_{j'}$ and
 $b^J_{i}$ contain no pair $(x_k, y_k)$ for $k=1$, $\dots$, $N$ except that 
 $a^{* J}_{j}$ and $a^{J}_j$ contain a pair $(x_j, y_j)$. Thus, by direct easy calculation,
 we obtain OPEs $a^{J}_{j'}(z) a^{* J}_{j}(w) \sim \hbar \delta_{j j'} / (z-w)$,
 $b^J_i(z) b^J_{i'}(w) \sim \hbar^2 \langle c_i, c_{i'} \rangle / (z-w)^2$ and all other
 combinations have trivial OPEs. Thus, we have the isomorphism of $\hbar$-adic vertex
 algebras of the statement.
\end{proof}
\end{proposition}

\section{Equivariant torus action and vertex algebra of global sections}
\label{sec:F-action}

In the previous sections, we defined the sheaf of $\hbar$-adic vertex algebras 
$\tldcalD^{ch}_{X, \hbar}$ over the hypertoric 
variety $\tldX$, and studied its structure. The space of global sections, 
$\tldcalD^{ch}_{X, \hbar}(\tldX)$
is naturally equipped with the structure of $\hbar$-adic vertex algebra. We also 
have an $\hbar$-adic vertex algebra constructed by the global BRST reduction
$H^{\infty/2+0}(\frg, \tldcalD^{ch}_{\frX, \hbar}(\tldfrX))$.
In this section, we construct the vertex algebras from these $\hbar$-adic vertex algebras using
a certain equivariant torus action, which reflect the essential structure of the original
$\hbar$-adic vertex algebras.

Consider an action of one-dimensional torus $\bbS = \C^{\times}$ on $\tldfrX$ which induces
an action on the structure sheaf $\calO_{\tldfrX} = \calO_{\frX} \otimes_{\C} \calO_{\frg^*}$
such that the weights of the generators with respect to the action is given by 
$\Swt(x_k) = \Swt(y_k) = 1/2$, $\Swt(c_i) = 1$ for $k=1$, $\dots$, $N$ and $i=1$, $\dots$, $M$.
Note that, with respect to this action, the Poisson bracket on $\calO_{\tldfrX}$ is homogeneous
of weight $-1$. Since the $\bbS$-action commutes with the $G$-action, we have the induced 
$\bbS$-action on the hypertoric variety $\tldX$.

Moreover, we have the equivariant $\bbS$-action on the sheaf $\tldcalD^{ch}_{\frX, \hbar}$
over $\C$ such that the weights of the generators given by 
$\Swt(x_{k (-n)}) = \Swt(y_{k (-n)}) = 1/2$, $\Swt(c_{i (-n)}) = 1$, $\Swt(\hbar) = 1$ 
and $\Swt(\bfone) = 0$
for $k=1$, $\dots$, $N$, $i=1$, $\dots$, $M$ and $n \in \Z_{\ge 1}$. Note that the OPEs of 
$\tldcalD^{ch}_{\frX, \hbar}$ are homogeneous with respect to the $\bbS$-action.
Extend this action onto the BRST complex $C_{\VA}$ by $\Swt(\psi^*_{i (-n)}) = 0$, 
$\Swt(\psi_{i (-n)}) = 1$ for $i=1$, $\dots$, $M$ and $n \in \Z_{\ge 1}$.
Then, the element $Q_{\VA} \in C_{\VA}$ is homogeneous of weight $\Swt(Q_{\VA}) = 1$, and
hence the coboundary operator $d_{\VA} = \hbar^{-1} Q_{\VA (0)}$ is a homogeneous 
operator of weight $0$ on the complex $C_{\VA}$. This implies that the BRST cohomology
sheaf $\calH^{\infty/2+\bullet}_{\VA}(\frg, \tldcalD^{ch}_{\frX, \hbar})$ is also equipped
with the induced equivariant $\bbS$-action over $\tldX$. In particular, the space of
global sections $\calH_{\VA}^{\infty/2+\bullet}(\frg, \tldcalD^{ch}_{\frX, \hbar})(\tldX)$ is
a $\C[[\hbar]]$-module with an $\bbS$-action over $\C$.

Recall the affine open covering $\tldfrX = \bigcup_{J} \tldfrU_J$. For any $J$, the open
subset $\tldfrU_J$ is closed under the $\bbS$-action, and $C_{\VA}(\tldfrU_J)$ is 
decomposed into the direct product of weight spaces because the coordinate functions of 
\refeq{eq:tld-local-trivial} are all homogeneous. Since the coboundary operators $d_{\VA}$
is homogeneous of weight $0$, the $0$-th cohomology group 
$\tldcalD^{ch}_{X, \hbar}(\tldU_J) = H^{0}(C_{\VA}(\tldfrU_J), d_{\VA})$
is also a direct product of weight spaces. Therefore, the $\hbar$-adic vertex algebra of global
sections can be decomposed into a direct product of weight spaces:
$\tldcalD^{ch}_{X, \hbar}(\tldX) = \prod_{m \ge 0} \tldcalD^{ch}_{X, \hbar}(\tldX)^{\bbS, m}$.
Note that the weights $m \in \frac{1}{2} \Z_{\ge 0}$ are non-negative and we have
$\tldcalD^{ch}_{X, \hbar}(\tldX)^{\bbS, 0} = \C \bfone$.
Consider the subspace
$\tldcalD^{ch}_{X, \hbar}(\tldX)_{fin} = \bigoplus_{m \in \frac{1}{2} \Z_{\ge 0}} \tldcalD^{ch}_{X, \hbar}(\tldX)^{\bbS, m}$.
This subspace is a $\C[\hbar]$-module since the weights are non-negative and 
$\Swt(\hbar)=1$. Moreover, 
since the OPEs preserve the $\bbS$-weight, they also preserve the subspace.
Now we set
\begin{equation}
 \label{eq:sfD}
  \tldsfD^{ch}(\tldX) =
\tldcalD^{ch}_{X, \hbar}(\tldX)_{fin} \bigm/ (\hbar - 1),
\end{equation}
the quotient space by the ideal generated by $\hbar - 1$.
It is a $\C$-vector space equipped with OPEs induced from ones on 
$\tldcalD^{ch}_{X, \hbar}(\tldX)$. Since the all 
identities between the vertex operators of 
$\tldcalD^{ch}_{X, \hbar}(\tldX)$ are satisfied
by the vertex operators of $\tldsfD^{ch}(\tldX)$, the $\C$-vector space $\tldsfD^{ch}(\tldX)$
is a vertex algebra.

Similarly, considering the $\hbar$-adic vertex subalgebra 
$H_{\VA}^{\infty/2 + 0}(\frg, \tldcalD^{ch}_{\frX, \hbar}(\tldfrX)) \subset \tldcalD^{ch}_{X, \hbar}(\tldX)$,
we have a $\C[\hbar]$-submodule
\[
 H_{\VA}^{\infty/2 + 0}(\frg, \tldcalD^{ch}_{\frX, \hbar}(\tldfrX))_{fin}
 = \bigoplus_{m \ge 0} H_{\VA}^{\infty/2 + 0}(\frg, \tldcalD^{ch}_{\frX, \hbar}(\tldfrX))^{\bbS, m}
 \subset H_{\VA}^{\infty/2 + 0}(\frg, \tldcalD^{ch}_{\frX, \hbar}(\tldfrX)).
\]
We define a vertex algebra over $\C$ by
\begin{equation}
 \label{eq:sfD-global}
  \sfD^{ch}(\tldX) = H_{\VA}^{\infty/2 + 0}(\frg, \tldcalD^{ch}_{\frX, \hbar}(\tldfrX))_{fin} 
  \bigm/ (\hbar - 1)
\end{equation}

\begin{definition}
 \label{def:hypertoric-VA}
 We call the vertex algebras $\tldsfD^{ch}(\tldX)$, $\sfD^{ch}(\tldX)$ defined by \refeq{eq:sfD}, 
 \refeq{eq:sfD-global}
 hypertoric vertex algebras. 
\end{definition}

\begin{remark}
 Later in \refprop{prop:global-VA}, we prove that the two vertex algebras
 $\tldsfD^{ch}(\tldX)$ and $\sfD^{ch}(\tldX)$ coincide.
\end{remark}

By the result of the previous section, the sheaf of $\hbar$-adic vertex algebra 
$\tldcalD^{ch}_{X, \hbar}$ is isomorphic to the tensor product of a $\beta\gamma$-system
and a Heisenberg vertex algebra. It gives an analog of Wakimoto realization
(free field realization) of the hypertoric
vertex algebra $\sfD^{ch}(\tldX)$ (and $\tldsfD^{ch}(\tldX)$). 
(cf. \cite{Wakimoto}, \cite{FF-Wakimoto})

For the affine open subset $\tldU_J \subset \tldX$, we have the restriction homomorphism
$\tldcalD^{ch}_{X, \hbar}(\tldX) \longrightarrow \tldcalD^{ch}_{X, \hbar}(\tldU_J)$ between
$\hbar$-adic vertex algebras. By \refprop{prop:local-descr}, we have
\[
 \tldcalD^{ch}_{X, \hbar}(\tldU_J) =
\C[[\hbar]][\,a^{* J}_{j (-n)}, a^{J}_{j (-n)}, b^{J}_{i (-n)} \;|\; j \not\in J, i=1, \dots, M, 
 n \in \Z_{\ge 1}\,]
\]
Then, the image of the $\C[\hbar]$-submodule $\calD^{ch}_{X, \hbar}$ under the homomorphism is
included in the $\C[\hbar]$-submodule 
$\C[\hbar][a^{* J}_{j (-n)}, a^{J}_{j (-n)}, b^{J}_{i (-n)} | i, j, n]$. 
Thus, we have the following $\C$-linear map
\begin{multline*}
 \sfD^{ch}(\tldX) \rightarrow \tldsfD^{ch}(\tldX) \rightarrow \\
 \C\bigl[\, a^{J*}_{j (-n)}, a^{J}_{j (-n)} \;\bigm\vert\; \substack{j \not\in J \\ n \in \Z_{\ge 1}}\,\bigr] 
 \otimes_{\C} \C\bigl[\,b_{i (-n)} \;\bigm\vert\; \substack{i=1, \dots, M \\ n \in \Z_{\ge 1}}\,\bigr] \bfone \\
 \simeq \calD^{ch}(\C^{(N-M)}) \otimes_{\C} V_{\langle, \rangle}(\frg).
\end{multline*}
by taking quotients by $(\hbar - 1)$ where $\calD^{ch}(\C^{(N-M)})$ is a 
$\beta\gamma$-system and $V_{\langle, \rangle}(\frg)$ is a Heisenberg vertex algebra. 
Clearly, this is a homomorphism between vertex algebras over $\C$.

For $\lambda \in \frg^*$, let $\pi_{\lambda}$ is the Heisenberg Fock space of
highest weight $\lambda$; i.e.  $\pi_{\lambda}$ is an irreducible highest weight 
module with a highest weight vector $| \lambda \rangle \in \pi_{\lambda}$ 
on which the action is given by 
$b^J_{i (0)} | \lambda \rangle = \lambda(c_i) | \lambda \rangle$ and
$b^J_{i (n)} | \lambda \rangle = 0$ for $i=1$, $\dots$, $M$ and $n > 0$.

\begin{proposition}
 \label{prop:Wakimoto}
 For each $J$ and $\lambda \in \frg^*$, we have an action of
 the hypertoric vertex algebra $\sfD^{ch}(\tldX)$ on 
 $\calF_{\beta\gamma} \otimes_{\C} \pi_{\lambda}$
 where $\calF_{\beta\gamma} = \C\bigl[\, a^{J*}_{j (-n)}, a^{J}_{j (-n)} \;\bigm\vert\; \substack{j \not\in J \\ n \in \Z_{\ge 1}}\,\bigr]$
 is the Fock space of the $\beta\gamma$-system and $\pi_{\lambda}$ is the Fock space
 of the Heisenberg vertex algebra of highest weight $\lambda$.
\end{proposition}

\section{Conformal vectors}
\label{sec:conf-vect}

In this section, we construct the conformal vector explicitly by an analog of the 
Segal-Sugawara construction.

First assume that the symmetric bilinear form $\langle \blkbar, \blkbar\rangle$ on $\frg$
is degenerate. In such a case, we have an element $\zeta = \sum_{i=1}^M a_i c_i \in \frg$ 
($a_i \in \C$ for $i=1$, $\dots$, $M$) such that $\langle \zeta, c_i \rangle = 0$ for any 
$i = 0$, $\dots$, $M$. Then, $\zeta = \zeta_{(-1)} \bfone \in C_{\VA}^0(\tldfrX)$ satisfies 
$\zeta_{(n)} \alpha = 0$ for any $\alpha \in C_{\VA}$ and $n \ge 0$. In particular, we have 
$d_{\VA} \zeta = \hbar^{-1} Q_{\VA(0)} \zeta = 0$. Clearly, 
$\zeta$ does not lie in $\Im d_{\VA}$ and thus $\zeta$ defines
a nonzero central vector in $H^0(C_{\VA}(\tldfrX), d_{\VA})$
and in $\sfD^{ch}(\tldX)$. Therefore, the vertex algebra $\sfD^{ch}(\tldX)$ has nontrivial
center and hence it is not a vertex operator algebra.

Now, assume that the symmetric bilinear form $\langle \blkbar, \blkbar \rangle$ is 
nondegenerate. Let $\{c^i\}_{i=1, \dots, M} \subset \frg$ be the dual basis of the basis 
$\{c_i\}_{i=1, \dots, M}$ with respect to the bilinear form. Set 
$\omega_H = (1/2) \sum_{i=1}^{M} c_{i (-1)} c^i \in V_{\frg, \langle, \rangle, \hbar}(\frg^*) \subset C^0_{\VA}(\tldfrX)$.
The following lemma is the standard fact.

\begin{lemma}
 \label{lemma:conf-vect-Heis}
 For $i=1$, $\dots$, $M$ and $m$, $n \in \Z$, we have 
 $[\omega_{H (m+1)}, c_{i (n)}] = - \hbar^2 n c_{i (m+n)}$. 
 In particular, one obtain $\mu_{ch}(A_i)_{(n)} \omega_H = - \hbar^2 \delta_{n1} c_i$ for
 $i=1$, $\dots$, $M$ and $n \ge 0$.
\end{lemma}

\begin{lemma}
 \label{lemma:Vir-Heis}
 We have the OPE
\[
 \omega_H(z) \omega_H(w) \sim 
 \frac{\hbar^4}{(z-w)^4} \frac{M}{2} + \frac{\hbar^2}{(z-w)^2} 2 \omega_H(w) +
 \frac{\hbar^2}{z-w} \partial_w \omega_H(w).
\]
\begin{proof}
 It is direct and standard calculation using \reflemma{lemma:conf-vect-Heis}.
\end{proof}
\end{lemma}

Let $\kappa \in \C$ be a parameter. For $j=1$, $\dots$, $N$ and $\kappa$, let
\[
 \omega_{\kappa, j} = \kappa x_{j (-2)} y_{j (-1)} \bfone + (\kappa - 1) x_{j (-1)} y_{j (-2)} \bfone
\in \calD^{ch}_{\frX, \hbar}(\frX) \subset C^0_{\VA}(\tldfrX).
\]

\begin{lemma}
 \label{lemma:conf-betagam}
 For $j$ $k=1$, $\dots$, $N$, and $m$, $n \in \Z$, we have
\[
  [\omega_{\kappa, j (m+1)}, x_{k (n)}] = - \hbar (n+\kappa) x_{j (m+n)}, \quad
 [\omega_{\kappa, j (m+1)}, y_{k (n)}] = - \hbar (n-\kappa+1) y_{j (m+n)}.
\]
In particular, we have
\[
 \mu_{ch}(A_i)_{(n)} \omega_{\kappa, j} =
 \begin{cases}
  \hbar \Delta_{ij} x_{j (-1)} y_j & (n=1) \\
  \hbar^2 \Delta_{ij} (1-2\kappa) \bfone & (n=2) \\
  0 & (\text{otherwise.})
 \end{cases}
\]
for $i=1$, $\dots$, $M$ and $n \ge 0$.
\end{lemma}

\begin{lemma}
 \label{lemma:Vir-betagamma}
 We have the OPE
\[
 \omega_{\kappa, j}(z) \omega_{\kappa, j}(w) \sim 
 \frac{\hbar^4}{(z-w)^4} \frac{-1}{2} + \frac{\hbar^2}{(z-w)^2} 2 \omega_{\kappa, j}(w) +
 \frac{\hbar^2}{z-w} \partial_w \omega_{\kappa, j}(w)
\]
for $k=1$, $\dots$, $N$ and for any $\kappa \in \C$.
\end{lemma}

Set $\omega_F = \sum_{i=1}^{M} \psi_{i (-2)}^* \psi_{i} \in Cl^{vert}(\frg \oplus \frg^*) \subset C^0_{\VA}(\tldfrX)$.
By direct calculation, we have the following lemma.

\begin{lemma}
 \label{lemma:conf-Clifford}
 We have the commutation relations $[\omega_{F (m+1)}, \psi^*_{i (n)}] = \hbar n \psi^*_{i (m+n)}$, 
 $[\omega_{F (m+1)}, \psi_{i (n)}] = \hbar n \psi_{i (m+n)}$ for $i=1$, $\dots$, $M$
 and $m$, $n \in \Z$. In particular, we have 
$d_{\VA} \omega_F = \sum_{i=1}^{M} \mu_{ch}(A_i)_{(-1)} \psi^*_{i (-2)} \bfone$.
Moreover, we have the following OPE
 \[
  \omega_F(z) \omega_F(w) \sim
 \frac{\hbar^4}{(z-w)^4} \frac{-2M}{2} + \frac{\hbar^2}{(z-w)^2} 2 \omega_F(w) +
 \frac{\hbar^2}{z-w} \partial_w \omega_F(w).
 \]
\end{lemma}

Now we set $\omega = \hbar \sum_{k=1}^{N} \omega_{1/2, k} + \omega_H + \hbar \omega_F \in C^0_{\VA}(\tldfrX)$.
Then the following proposition is obvious from the above lemmas.

\begin{proposition}
 \label{prop:conf-vect-def}
 We have $d_{\VA}(\omega) = 0$, and thus $\omega \in C_{\VA}^0(\tldfrX)$ defines an element in 
 $H^0(C_{\VA}(\tldfrX), d_{\VA})$ and in
 $\sfD^{ch}(\tldfrX)$ which we also write $\omega$. Moreover, the element $\omega$ has the OPE
 \[
   \omega(z) \omega(w) \sim 
 \frac{\hbar^4}{(z-w)^4} \frac{-M-N}{2} + \frac{\hbar^2}{(z-w)^2} 2 \omega(w) +
 \frac{\hbar^2}{z-w} \partial_w \omega(w).
 \]
 Namely, $\omega \in \sfD^{ch}(\tldX)$ is a conformal vector.
\end{proposition}

The operator $\omega_{(1)}$ gives a non-negative grading on $C_{\VA}(\tldfrX)$;
$\confwt(x_k) = \confwt(y_k) = 1/2$, $\confwt(c_i) = 1$, $\confwt(\psi^*_i) = 0$
and $\confwt(\psi_i) = 1$ for $k=1$, $\dots$, $N$, $i=1$, $\dots$, $M$. 
The vertex algebra $\sfD^{ch}(\tldX)$ is $\frac{1}{2}\Z_{\ge 0}$-graded by the action of $\omega_{(1)}$ 
such that any element of conformal weight $0$ is proportional to the vacuum $\bfone$.
Therefore, $\sfD^{ch}(\tldX)$ is a vertex operator algebra. 

Moreover, take $\lambda = (\lambda_k)_{k=1, \dots, N} \in \R^N$, an orthogonal vector with 
all row vectors $\Delta^i \defeq (\Delta_{ij})_{j=1, \dots, N}$ of the matrix $\Delta$ for $i = 1$, $\dots$, $M$.
Then, the vector 
\[
 \omega_{\lambda} = \omega + \hbar \sum_{k=1}^{N} 
 \lambda_k (x_{k (-2)} y_{k (-1)} + x_{k (-1)} y_{k (-2)}) \bfone
 = \hbar \sum_{k=1}^{N} \omega_{(1/2+\lambda_k), k} + \omega_H + \hbar \omega_F
\]
is also a conformal vector in $C_{\VA}(\tldfrX)$. Since $\lambda$ is orthogonal with 
$\Delta^i = (\Delta_{ij})_{j=1, \dots, N}$ for $i=1$, $\dots$, $M$, we have
\[
 \mu_{ch}(A_i)_{(n)} \sum_{k=1}^{N}  
\lambda_k (x_{k (-2)} y_{k (-1)} + x_{k (-1)} y_{k (-2)}) \bfone
= \delta_{n 2} \hbar^2 \sum_{k=1}^{N} \Delta_{ik} (-2 \lambda_k) \bfone = 0
\]
for all $i=1$, $\dots$, $N$ and $n \ge 0$. Thus, we have $\omega_{\lambda} \in \Ker d_{\VA}$,
and hence $\omega_{\lambda}$ induces a conformal vector in $\sfD^{ch}(\tldX)$.
Note that $\omega_{\lambda (1)}$ also gives a grading on $\sfD^{ch}(\tldX)$ but the 
grading may be negative in contrast to the standard one $\omega_{(1)}$.

\section{Zhu algebras}
\label{sec:Zhu-algebra}

In this section, we discuss the Zhu algebra of the hypertoric vertex algebra,
an associative algebra which reflects fundamental aspects of the representation theory of
the corresponding vertex operator algebra. Our goal is to show that the Zhu algebra of
$\sfD^{ch}(\tldX)$ coincides with the universal family of quantization of the Poisson 
algebra $\C[X]$.

\subsection{The definition of Zhu algebras}
\label{sec:define-Zhu}

Let $V = \bigoplus_{\Deg \ge 0} V_\Deg$ be a vertex algebra with $\Z_{\ge 0}$-grading. For 
a homogeneous element $a \in V_\Deg$, we denote its grading $\Deg_a = \Deg$. For a 
homogeneous element $a \in V_{\Deg_a}$, an element $b \in V$ and positive integer 
$m \in \Z_{> 0}$, we define
\[
 a *_m b = \sum_{j \ge 0} \binom{\Deg_a}{j} a_{(-m+j)} b \in V,
\]
and extend it on $V$ linearly.
We simply denote $* = *_{1}$, $\circ = *_{2}$ for $m=1$, $2$. Let 
$A(V) = V \,/\, V \circ V$ be the quotient vector space where 
$V \circ V = \{ a \circ b \;|\; a, b \in V \}$. As proved in \cite{Zhu}, 
the product $* = *_1$ is a linear associative product on
the vector space $A(V)$ with the unit $\bfone$. The associative algebra 
$A(V)$ is called the Zhu algebra of the vertex algebra $V$.

Besides the Zhu algebra $A(V)$, we also have a Poisson algebra corresponding to
the vertex algebra $V$. Consider the vector space $\Abar(V) = V \,/\, V_{(-2)} V$
where $V_{(-2)} V = \{ a_{(-2)} b \;|\; a, b \in V \}$. The vertex algebra operator ${}_(-1)$
gives a commutative associative product on $\Abar(V)$, and moreover, $\Abar(V)$ is
a Poisson algebra with the Poisson bracket $\{a, b\} = a_{(0)} b$ modulo $V_{(-2)} V$.
We call the Poisson algebra $\Abar(A)$ the $C_2$ Poisson algebra of the vertex algebra $V$.
Note that, while the definition of Zhu algebra $A(V)$ requires the $\Z_{\ge 0}$-grading
on the vertex algebra $V$, the grading is not needed to define  $C_2$ Poisson algebra $\Abar(V)$.
In some known cases, the Zhu algebra gives a quantization of the $C_2$ Poisson algebra;
e.g. the affine vertex operator algebra associated with the simple Lie algebra, 
Virasoro vertex algebra and $\beta\gamma$ systems. 

Also for an $\hbar$-adic vertex algebra $V_{\hbar}$, we define 
$\Abar(V_{\hbar}) = V_{\hbar} \,/\, V_{\hbar (-2)} V_{\hbar}$, a commutative 
$\C[[\hbar]]$-algebra. For the sheaf of $\hbar$-adic vertex algebras 
$\tldcalD^{ch}_{X, \hbar}$ over $\tldX$, we define the sheaf of $\C[[\hbar]]$-algebras
$\calAbar(\tldcalD^{ch}_{X,\hbar})$ as the quotient sheaf 
$\calAbar(\tldcalD^{ch}_{X, \hbar}) = \tldcalD^{ch}_{X, \hbar} \,/\, \tldcalD^{ch}_{X, \hbar (-2)} \tldcalD^{ch}_{X, \hbar}$.
Namely, it is the sheaf associated with the presheaf 
$\tldU \mapsto \Abar(\tldcalD^{ch}_{X, \hbar}(\tldU))$ for an open subset $\tldU \subset \tldX$.

\subsection{Quantization of the hypertoric variety}
\label{sec:quant-hypertoric}

The associative algebra quantizing the hypertoric variety $X$ was first introduced by
I.~Musson and M.~Van den Bergh in \cite{Musson-Van-der-Burgh}. 

Let $\calD(V)$ be the Weyl algebra on the affine space $V = \C^N$, that is the algebra
of differential operators with polynomial coefficients. We denote the standard coordinate
functions on $V$ by $x_1$, $\dots$, $x_N$ as in \refsec{sec:hypertoric-varieties},
and the corresponding differential operators $\partial_k = \partial / \partial x_k$
for $k=1$, $\dots$, $N$. Then the Weyl algebra $\calD(V)$ is isomorphic to
$\C[x_k, \partial_k \;|\; k=1, \dots, N]$ as a $\C$-vector space.
The action of the torus $G = (\C^\times)^M$ on $V$ induces an action on the algebra $\calD(V)$.
Define a map $\mu_D : \frg \longrightarrow \calD(V)$ by 
$A_i \mapsto \mu_D(A_i) = \sum_{k=1}^{N} \Delta_{ik} x_k \partial_k$ for $i=1$, $\dots$, $M$.
Clearly, this map quantizes the comoment map $\mu^*$ and we call $\mu_D$ a quantized comoment
map. Set $\tldcalD(V) = \calD(V) \otimes_{\C} \C[c_1, \dots, c_M]$, and extend 
the action of the torus $G$ onto $\tldcalD(V)$ so that $G$ acts on $\C[c_1, \dots, c_M]$
trivially. Define the associative algebra $\sfD(\tldX)$ by quantum Hamiltonian reduction
as follows:
\begin{equation}
\label{eq:def-q-hypertoric} 
 \sfD(\tldX) = \Bigl(\calD(V) \bigm/ \sum_{i=1}^{M} \calD(V) (\mu_D(A_i) - c_i) \Bigr)^G
 = \calD(V)^G \bigm/ \sum_{i=1}^{M} \calD(V)^G (\mu_D(A_i) - c_i).
\end{equation}
It is not difficult to examine that $\sfD(\tldX)$ is an associative algebra,
and its associated graded algebra with respect 
to the Berenstein filtration, i.e. the filtration induced from 
$\deg x_k = \deg \partial_k = 1$ and  $\deg c_i = 0$, coincides with $\C[\tldX]$ 
as Poisson algebras. The algebra $\sfD(\tldX)$ is an algebra over 
$\C^M = \Spec \C[c_1, \dots, c_M]$, and it is a family of filtered quantizations of
the Poisson algebra $\C[X]$, while the Poisson algebra $\C[\tldX]$ is a family of
Poisson deformations of $\C[X]$ over $\C^M$ in the sense of \cite{Losev}, \cite{Losev-2016}.

The algebra $\sfD(\tldX)$ was introduced in \cite{Musson-Van-der-Burgh},
and it is called a quantized hypertoric algebra or a hypertoric enveloping algebra.
One can construct a sheaf of associative $\C[[\hbar]]$-algebras on $\tldX$ whose algebra
of global sections coincides with $\sfD(\tldX)$. See \cite{Bellamy-K} and \cite{BLPW}.
Moreover we can describe the above quantum Hamiltonian reduction by a certain BRST 
cohomology, which is analogous to the BRST cohomology in this paper. See \cite{BRST}.

Consider the action of the $N$-dimensional abelian Lie algebra 
$\C^N = \bigoplus_{k=1}^{N} \C x_k \partial_k$ on $\calD(V)$ by the commutation
$[x_k \partial_k, \blkbar]$ for $k=1$, $\dots$, $N$. The action corresponds an
action of the $N$-dimensional torus $\bbT = (\C^{\times})^N$ on $\calD(V)$
induced from the natural action on $\C^N$. The algebra $\calD(V)$ is decomposed into
the direct sum of weight spaces with respect to this action: 
$\calD(V) = \bigoplus_{\zeta \in \Z^N} \calD(V)^{\bbT, \zeta}$. 
Consider the sublattice $\bigoplus_{i=1}^{M} \Z \Delta^i \subset \Z^N$ where
$\Delta^i = (\Delta_{ij})_{j=1, \dots, N}$. It can be identified with the weight lattice of 
the torus $G$ and its Lie algebra $\frg = \bigoplus_{i=1}^M \C A_i$
because $A_i \in \frg$ acts on $\calD(V)$ by 
$\mu_D(A_i) = \sum_{j=1}^N \Delta_{ij} x_j \partial_j$.
Take the orthogonal sublattice $\Lambda_0 \subset \Z^N$ of $\bigoplus_{i=1}^{M} \Z \Delta^i$.
Then, we have $\calD(V)^G = \bigoplus_{\zeta \in \Lambda_0} \calD(V)^{\bbT, \zeta}$ and it
induces the weight decomposition of the quantized hypertoric algebra: 
$\sfD(\tldX) = \bigoplus_{\zeta \in \Lambda_0} \sfD(\tldX)^{\bbT, \zeta}$. The following
lemma is obvious.

\begin{lemma}
 \label{lemma:4}
 The weight space $\sfD(\tldX)^{\bbT} \subset \sfD(\tldX)$ of weight $0$ is given by
 \[
  \sfD(\tldX)^{\bbT} = \C[x_1 \partial_1, \dots, x_N \partial_N] \otimes_{\C} \C[c_1, \dots, c_M].
 \]
 Setting
 $P_{\zeta} = \prod_{k: \zeta_k > 0} x_k^{\zeta_k} \prod_{k: \zeta_k < 0} \partial_k^{-\zeta_k}$
 for $\zeta = (\zeta_1, \dots, \zeta_N) \in \Lambda_0$, the weight space 
 $\sfD(\tldX)^{\bbT, \zeta}$ is a $\sfD(\tldX)^{\bbT}$-module generated by $P_{\zeta}$.
\end{lemma}

Clearly, the associated graded algebra $\C[\tldX]$ has also the same weight decomposition:
$\C[\tldX] = \bigoplus_{\zeta \in \Lambda_0} \C[\tldX]^{\bbT, \zeta}$. For each $\zeta \in \Lambda_0$,
we have the same description for $\C[\tldX]^{\bbT, \zeta}$ as \reflemma{lemma:4}; that is,
$\C[\tldX]^{\bbT} = \C[x_1 y_1, \dots, x_N y_N, c_1, \dots, c_M]$ for $\zeta = 0$ and 
$\C[\tldX]^{\bbT, \zeta}$ is a $\C[\tldX]^{\bbT}$-module generated by $P_{\zeta}$, where we identify 
$P_{\zeta} \in \sfD(\tldX)^{\bbT, \zeta}$ with its image 
$P_{\zeta} = \prod_{k: \zeta_k > 0} x_k^{\zeta_k} \prod_{k: \zeta_k < 0} y_k^{-\zeta_k} \in \C[\tldX]^{\bbT, \zeta}$.

\subsection{Weyl group symmetries}
\label{sec:Weyl-group}

The Weyl algebra $\calD(V)$ has natural automorphisms in $(\Z/2\Z)^N \ltimes \frS_N$, 
generated by permutations $\sigma \in \frS_N$, $\sigma(x_k) = x_{\sigma(k)}$, 
$\sigma(\partial_k) = \partial_{\sigma(k)}$, and Fourier transformations
$x_k \mapsto - y_k$, $y_k \mapsto x_k$, for each $k=1$, $\dots$, $N$. It naturally
induces an action on the weight lattice $\Z^N$. 
Let $\bbW$ be the subgroup of all elements in $(\Z/2\Z)^N \ltimes \frS_N$ which fix 
the sublattice $\Lambda_0$ pointwise. Since $\Delta^1$, $\dots$, $\Delta^M$ span
the sublattice which is orthogonal to $\Lambda_0$, an element $\sigma \in \bbW$ maps
$\mu_D(A_i)$ to a linear combination $\sum_{j=1}^{M} \lambda_j \mu_D(A_j)$, 
$\lambda_j \in \Z$. 
Then, the action of $\bbW$ on $\calD(V)$ is extended
onto $\tldcalD(V) = \calD(V) \otimes_{\C} \C[c_1, \dots, c_M]$ by 
$\sigma(c_i) = \sum_{j=1}^{M} \lambda_j c_j$. By the definition \refeq{eq:def-q-hypertoric},
the action $\bbW$ on $\tldcalD(V)$ induces automorphisms of the quantized hypertoric algebra 
$\sfD(\tldX)$. It also induces automorphisms of the Poisson algebra $\C[\tldX]$.
The algebras $\sfD(\tldX)$ and $\C[\tldX]$ also have other automorphisms which fix the parameters
$c_1$, $\dots$, $c_M$, denoted $\bbV$ in \cite[Section 8.1]{BLPW}, but we will ignore
such automorphisms. Now consider the $\bbW$-invariant subalgebras $\sfD(\tldX)^{\bbW}$
and $\C[\tldX]^{\bbW}$. The algebra $\sfD(\tldX)^{\bbW}$ (resp. $\C[\tldX]^{\bbW}$) is also a 
family of filtered quantizations (resp. Poisson deformations) of the Poisson algebra $\C[X]$
over the space $\C^M / \bbW$. By Corollary~2.13 and Proposition~3.5 in \cite{Losev-2016},
$\sfD(\tldX)^{\bbW}$ (resp. $\C[\tldX]^{\bbW}$) is characterized as the universal family of 
filtered quantizations (resp. Poisson deformations) of the Poisson algebra $\C[X]$.

Using \reflemma{lemma:4} we have description of the $\bbW$-invariant subalgebra 
$\sfD(\tldX)^{\bbW}$ as follows: 
By the orthogonal decomposition $\bigoplus_{i=1}^{M} \Z \Delta^i \oplus \Lambda_0$,
for $k=1$, $\dots$, $N$, 
we have the decomposition $x_k \partial_k = \sum_{i=1}^{M} \beta_i \mu_D(A_i) + z$ 
where $\beta_i \in \C$ and $z \in \bigoplus_{k=1}^{N} \C x_k \partial_k$ is an element which is orthogonal to
$\mu_D(A_i)$ for all $i=1$, $\dots$, $M$. Set 
\begin{equation}
 \label{eq:euler-op}
 H_k = x_k \partial_k - \sum_{i=1}^{M} \beta_i c_i \in \sfD(\tldX)
\end{equation}
for $k=1$, $\dots$, $N$. Since $H_k = z$ in $\sfD(\tldX)$ and the group $\bbW$ 
fixes $\Lambda_0$ pointwise, $H_k$ is invariant under the action of $\bbW$ on $\sfD(\tldX)$.
Next, consider the element $P_{\zeta} \in \sfD(\tldX)^{\bbT, \zeta}$ in \reflemma{lemma:4}. 
Since $\sigma \in \bbW$ fixes the sublattice $\Lambda_0$ pointwise, $\sigma(P_{\zeta})$ is
again an element of $\sfD(\tldX)^{\bbT, \zeta}$. Moreover, we have $\sigma(P_{\zeta}) = P_{\zeta}$
since $P_{\zeta}$ is the only element which has none of the factors $x_k \partial_k$ for
any $k=1$, $\dots$, $N$. Therefore, $P_{\zeta}$ is a $\bbW$-invariant element in 
$\sfD(\tldX)^{\bbT, \zeta}$.

\begin{lemma}
 \label{lemma:gen-of-W-inv}
 The set of polynomials $\{P_{\zeta}\,|\, \zeta \in \Lambda_0 \} \cup \{ H_k \,|\, k=1, \dots, N\}$
 generates the $\bbW$-invariant subalgebra $\C[\tldX]^{\bbW}$.
 \begin{proof}
  Let $R$ be a subalgebra of $\C[\tldX]$ generated by the elements 
  $\{P_{\zeta}\,|\, \zeta \in \Lambda_0 \} \cup \{ H_k \,|\, k=1, \dots, N\}$.
  Since the generators are $\bbW$-invariant and homogeneous, the subalgebra $R$ is 
  a graded subalgebra of $\C[\tldX]^{\bbW}$. Set $S = R \cap \C[\frg^*]^{\bbW}$. Then, we
  have $R \otimes_{S} \C \simeq \C[X]$ where $\C$ is an $S$-algebra induced from the
  specialization $c_i \mapsto 0$ for $i=1$, $\dots$, $M$. Thus, $R$ is a graded family of
  Poisson deformation of $\C[X]$ over $S$. By \cite[Proposition 2.12]{Losev-2016}, 
  we have a unique homomorphism $\C[\frg^*]^{\bbW} \longrightarrow S$ which induces an
  isomorphism $\C[\tldX]^{\bbW} \otimes_{\C[\frg^*]^{\bbW}} S \simeq R$ intertwining
  the isomorphisms 
  $R \otimes_{S} \C \isoto \C[X] \isoto \C[\tldX]^{\bbW} \otimes_{\C[\frg^*]^{\bbW}} \C$.
  By the definition of $R$, the embedding $R \hookrightarrow \C[\tldX]^{\bbW}$ also 
  intertwines the isomorphisms 
  $R \otimes_{S} \C \simeq \C[X] \simeq \C[\tldX]^{\bbW} \otimes_{\C[\frg^*]^{\bbW}} \C$.
  Consider the composition $\varphi: \C[\tldX]^{\bbW} \longrightarrow \C[\tldX]^{\bbW}$ 
  of the above homomorphisms $\C[\tldX]^{\bbW} \longrightarrow R$ and 
  $R \hookrightarrow \C[\tldX]^{\bbW}$. Then, $\varphi$ intertwines the isomorphisms
  $\C[\tldX]^{\bbW} \otimes_{\C[\frg^*]^{\bbW}} \C \simeq \C[X] \simeq \C[\tldX]^{\bbW} \otimes_{\C[\frg^*]^{\bbW}} \C$.
  Therefore, the homomorphism $\varphi$ is an isomorphism by the universality.
  This implies $R = \C[\tldX]^{\bbW}$.
 \end{proof}
\end{lemma}

\subsection{The $C_2$ Poisson algebra}
\label{sec:C_2-Poisson}

Now we determine the $C_2$ Poisson algebra $\Abar(\sfD^{ch}(\tldX))$ of the hypertoric
vertex algebra $\sfD^{ch}(\tldX)$. Consider the affine open covering 
$\tldX = \bigcup_{J} \tldU_J$, and we have an isomorphism of \refprop{prop:local-descr}:
\begin{align*}
 \tldcalD^{ch}_{X, \hbar}(\tldU_J) &= \C[[\hbar]]\bigl[\, a^{*}_{j (-n)}, a_{j (-n)} \,\bigm|\,
 \substack{ j \not\in J \\ n \ge 1} \,\bigr] \otimes_{\C[[\hbar]]}
 \C[[\hbar]]\bigl[\, b_{i (-n)} \,\bigm|\, \substack{i=1, \dots, M \\ n \ge 1} \,\bigr] \\
 & \simeq \calD^{ch}(\C^{2(N-M)})_{\hbar} \otimes_{\C[[\hbar]]} V_{\langle, \rangle, \hbar }(\frg).
\end{align*}
Thus, its $C_2$ Poisson algebra 
$\Abar(\tldcalD^{ch}_{X, \hbar}(\tldU_J)) = \calAbar(\tldcalD^{ch}_{X, \hbar})(\tldU_J)$ 
for each affine open subset $\tldU_J \subset \tldX$ is given by
\begin{multline*}
\calAbar(\tldcalD^{ch}_{X, \hbar})(\tldU_J) = \Abar(\tldcalD^{ch}_{X, \hbar}(\tldU_J))
 \\ 
\simeq \C[[\hbar]][ a^*_{j (-1)}, a_{j (-1)} \,|\, j \not\in J ] \otimes_{\C[[\hbar]]}
 \C[[\hbar]][ b_{i (-1)} \,|\, i=1, \dots, M ]
 \simeq \calO_{\tldX}(\tldU_J).
\end{multline*}
Moreover, the coordinate transformation of $\calD^{ch}_{X, \hbar}$ on $\tldU_I \cap \tldU_J$
maps $b_{i (-1)}^{I}$ to 
$b_{i (-1)}^{J} - \sum_{j=1}^{M} \langle c_i, c_j \rangle \partial \log(T^I_j / T^J_j)$
for $i=1$, $\dots$, $M$
and the local sections $\partial \log(T^I_j / T^J_j) \equiv 0$ in the $C_2$ Poisson algebra
$\Abar(\tldcalD^{ch}_{X, \hbar}(\tldU_I \cap \tldU_J))$. Thus, this coordinate transformation 
induces the coordinate transformation of $\calAbar(\tldcalD^{ch}_{X, \hbar})$ such that
$b^I_{i (-1)}$ is mapped to $b^J_{i (-1)}$ for each $i=1$, $\dots$, $M$ and each 
$\tldU_{I} \cap \tldU_{J}$.

\begin{lemma}
 \label{lemma:3}
 We have an isomorphism of sheaves of $\C[[\hbar]]$-algebras
$\calO_{\tldX}[[\hbar]] \longrightarrow \calAbar(\tldcalD^{ch}_{X, \hbar})$ which is
 locally given by
\begin{align*}
 \calO_{\tldX}(\tldU_J)[[\hbar]] &\longrightarrow \calAbar(\tldcalD^{ch}_{X, \hbar})(\tldU_J) & \\
 a^{J *}_{j},\, a^{J}_{j} &\mapsto a^{J*}_{j (-1)},\, a^{J}_{j (-1)}, & (j \not\in J) \\
 c_i &\mapsto b^{J}_{i (-1)}. & (i=1, \dots, M)
\end{align*}
\end{lemma}

Since the global section functor $\Gamma(\tldX, \blkbar)$ is left adjoint, 
$\Abar(\tldcalD^{ch}_{X, \hbar}(\tldX))$ is a subalgebra of 
$\calAbar(\tldcalD^{ch}_{X, \hbar})(\tldX) \simeq \calO_{\tldX}(\tldX)[[\hbar]]$. 
From this fact, we obtain the following fundamental fact for the hypertoric vertex algebra.

\begin{proposition}
 \label{prop:global-VA}
 We have $\tldcalD^{ch}_{X, \hbar}(\tldX) = H^{\infty/2+0}_{\VA}(\frg, \tldcalD^{ch}_{\frX, \hbar}(\tldfrX)) = H^0(C_{\VA}(\tldfrX), d_{\VA})$, 
 and hence $\tldsfD^{ch}(\tldX) = \sfD^{ch}(\tldX)$.
\begin{proof}
 If $H^0(C_{\VA}(\tldfrX), d_{\VA}) \ne \tldcalD^{ch}_{X, \hbar}(\tldX)$, then clearly
 there exists an element of the $C_2$ Poisson algebra $\Abar(\tldcalD^{ch}_{X, \hbar}(\tldX))$
 which does not lie in the image of 
 $H^0(C_{\VA}(\tldfrX), d_{\VA}) \subset \tldcalD^{ch}_{X, \hbar}(\tldX)$.
 However, $\Abar(\tldcalD^{ch}_{X, \hbar}(\tldX))$ is a subalgebra of 
 $\calO_{\tldX}(\tldX)[[\hbar]]$ and any element of 
 $\calO_{\tldX}(\tldX) \simeq \C[\tldfrX]^{G} / \sum_{i} \C[\tldfrX]^{G} (\mu^*(A_i) - c_i)$
 is represented by an element of $\C[\tldfrX]$. Thus, we have no element in 
 $\Abar(\tldcalD^{ch}_{X, \hbar}(\tldX))$ which does not lie in the image of
 $H^0(C_{\VA}(\tldfrX), d_{\VA})$.
\end{proof}
\end{proposition}

Recall the definition
$\sfD^{ch}(\tldX) = \tldsfD^{ch}(\tldX) = \tldcalD^{ch}_{X, \hbar}(\tldX)_{fin} / (\hbar - 1)$. 
By the isomorphism theorem, we have
\begin{multline}
 \label{eq:2}
 \Abar(\sfD^{ch}(\tldX)) = \Abar(\tldcalD^{ch}_{X, \hbar}(\tldX)_{fin} / (\hbar - 1)) \\
\simeq \Abar(\tldcalD^{ch}_{X, \hbar}(\tldX)_{fin})  / (\hbar - 1) 
\subset \calO_{\tldX}(\tldX)[\hbar] / (\hbar - 1) \simeq \C[\tldX].
\end{multline}

Now recall the element 
$P_{\zeta} = \prod_{k: \zeta_k > 0} x_k^{\zeta_k} \prod_{k: \zeta_k < 0} \partial_k^{-\zeta_k} \in \sfD(\tldX)^{\bbT, \zeta}$ 
for $\zeta \in \Lambda_0$ in \reflemma{lemma:4}. We consider the corresponding element 
\[
\tldP_{\zeta} = \prod_{k: \zeta_k > 0} x_{k (-1)}^{\zeta_k} \prod_{k: \zeta_k < 0} y_{k (-1)}^{- \zeta_k} \bfone \in C^0_{\VA}(\tldfrX) 
\]
of the BRST complex. Since $\zeta \in \Lambda_0$ is
orthogonal to $\Delta^i$ for $i=1$, $\dots$, $M$ and the element $\tldP_{\zeta}$ 
has none of the factors $x_{k (-1)} y_{k (-1)}$ for $k=1$, $\dots$, $N$, we have
$d_{\VA}(\tldP_{\zeta}) = 0$. Thus $\tldP_{\zeta}$ defines an element in 
$\sfD^{ch}(\tldX)$, and in its $C_2$ Poisson algebra 
$\Abar(\sfD^{ch}(\tldX))$. We denote these elements the same notation 
$\tldP_{\zeta}$. Next, recall the element $H_k$ for $k=1$, $\dots$, $N$ in \refeq{eq:euler-op}.
We define the corresponding element 
\[
 \tldH_k = x_{k (-1)} y_{k} - \sum_{i=1}^{M} \beta_i c_i \in C_{\VA}^0(\tldfrX)
\]
for $k=1$, $\dots$, $N$. Since $H_k \equiv z$ is orthogonal to $\mu_D(A_i)$ in 
$\bigoplus_{j=1}^N \C x_j \partial_j$ for all $i=1$, $\dots$, $M$, we have
$\mu_{ch}(A_i)_{(n)} \tldH_k = 0$ for all $n \ge 0$, and hence $d_{\VA}(\tldH_k) = 0$.
We denote the corresponding element in $\sfD^{ch}(\tldX)$ and 
$\Abar(\sfD^{ch}(\tldX))$ the same notation $\tldH_k$.
Clearly, $\tldH_1$, $\dots$, $\tldH_k$ together with the radical of the bilinear form 
$\langle \blkbar, \blkbar\rangle$ on $\bigoplus_{i=1}^{M} \C c_i \subset C_{\VA}^0(\tldfrX)$
span the image of the space 
$\bigoplus_{k=1}^{N} \C x_{k (-1)} y_k \oplus \bigoplus_{i=1}^{M} \C c_i$ in 
$\tldcalD^{ch}_{X, \hbar}(\tldX)$. By \refeq{eq:2} and \reflemma{lemma:gen-of-W-inv}, 
we have the following proposition.

\begin{proposition}
 \label{prop:C2-Poisson}
 The $C_2$ Poisson algebra $\Abar(\sfD^{ch}(\tldX))$ of the hypertoric vertex algebra 
 $\sfD^{ch}(\tldX)$ is a subalgebra of $\C[\tldX]$ which contains 
 the $\bbW$-invariant subalgebra $\C[\tldX]^{\bbW}$, 
 under the identification given by $\tldH_k \mapsto H_k$ 
 for $k=1$, $\dots$, $N$ and $\tldP_{\zeta} \mapsto P_{\zeta}$ for $\zeta \in \Lambda_0$.
\end{proposition}

\subsection{Zhu algebra}

As the final goal of the present paper, we determine the Zhu algebra 
$A(\sfD^{ch}(\tldX))$ of the hypertoric vertex algebra $\sfD^{ch}(\tldX)$.

Consider a $\frac{1}{2} \Z_{\ge 0}$-graded vertex algebra structure on the BRST
complex $C_{\VA}(\tldfrX)$, given by $\Deg_{x_k} = \Deg_{y_k} = 1/2$, 
$\Deg_{c_i} = 1$, $\Deg_{\psi^*_i} = 0$ and $\Deg_{\psi_i} = 1$ for $k=1$, $\dots$, $N$
and $i=1$, $\dots$, $M$. This grading is compatible with the conformal weights on
$C_{\VA}(\tldfrX)$ introduced in \refsec{sec:conf-vect} when the bilinear form 
$\langle \blkbar, \blkbar \rangle$ on $\bigoplus_{i=1}^{M} \C c_i$ is nondegenerate.
Thus, the coboundary operator $d_{\VA}$ is homogeneous of degree $0$, and hence
$\tldcalD^{ch}_{X, \hbar}(\tldX)$ and $\sfD^{ch}(\tldX)$ are also 
$\frac{1}{2} \Z_{\ge 0}$-graded. Using this grading, we define the Zhu algebra
$A(\sfD^{ch}(\tldX))$ of the hypertoric vertex algebra $\sfD^{ch}(\tldX)$.

First we characterize $A(\sfD^{ch}(\tldX))$ as a quantization of the $C_2$ Poisson
algebra $\Abar(\sfD^{ch}(\tldX))$.
Recall that the hypertoric vertex algebra 
$\sfD^{ch}(\tldX) = \tldcalD^{ch}_{X, \hbar}(\tldX)_{fin} / (\hbar - 1)$ is equipped
with a filtration induced from the $\hbar$-adic filtration on $\tldcalD^{ch}_{X, \hbar}$.
The filtration induces a filtration of the associative algebra $A(\sfD^{ch}(\tldX))$.

\begin{proposition}
 \label{prop:Zhu-quantization}
 The Zhu algebra $A(\sfD^{ch}(\tldX))$ is a quantization of the $C_2$ Poisson algebra
 $\Abar(\sfD^{ch}(\tldX))$. Namely, the associated graded algebra of 
 $A(\sfD^{ch}(\tldX))$ with respect to the above filtration 
 is isomorphic to $\Abar(\sfD^{ch}(\tldX))$ as a Poisson algebra over $\C$.
 \begin{proof}
  Note that the $\C[[t]]$-algebra $A(\tldcalD^{ch}_{X, \hbar}(\tldX))$ is the Rees algebra
  of the filtered algebra 
  $A(\sfD^{ch}(\tldX)) \simeq A(\tldcalD^{ch}_{X, \hbar}(\tldX)_{fin}) / (\hbar - 1)$.
  Thus, the associated graded algebra with respect to the filtration is given by
  $\Gr A(\sfD^{ch}(\tldX)) \simeq A(\tldcalD^{ch}_{X, \hbar}(\tldX)) / (\hbar)  \simeq A(\tldcalD^{ch}_{X, \hbar}(\tldX) / (\hbar))$. 
  In the commutative vertex algebra $\tldcalD^{ch}_{X, \hbar}(\tldX) / (\hbar)$, 
  we have $a \circ b = a_{(-2)} b + \Deg_a a_{(-1)} b$ and
  $a * b = a_{(-1)} b$ for $a$, 
  $b \in \tldcalD^{ch}_{X, \hbar}(\tldX) / (\hbar)$ where $d_a \in \R$ is the degree of $a$.
  Thus, the Zhu algebra $A(\tldcalD^{ch}_{X, \hbar}(\tldX) / (\hbar))$ is isomorphic
  to the $C_2$ Poisson algebra $\Abar(\tldcalD^{ch}_{X, \hbar}(\tldX) / (\hbar))$.
  By the isomorphism theorem, we have
\begin{multline*}
  \Abar(\tldcalD^{ch}_{X, \hbar}(\tldX) / (\hbar)) = \Abar(\tldcalD^{ch}_{X, \hbar}(\tldX)_{fin} / (\hbar)) 
 \simeq \Abar(\tldcalD^{ch}_{X, \hbar}(\tldX)_{fin}) / (\hbar) \\
  \simeq \Gr \Abar(\tldcalD^{ch}_{X, \hbar}(\tldX)_{fin}) / (\hbar - 1)
  \simeq \Gr \Abar(\sfD^{ch}(\tldX)) \simeq \Abar(\sfD^{ch}(\tldX)).
\end{multline*}
 \end{proof}
\end{proposition}

By \cite[Proposition 3.5]{Losev-2016}, the $\bbW$-invariant subalgebra $\sfD(\tldX)^{\bbW}$
of the quantized hypertoric algebra $\sfD(\tldX)$ gives a universal family of
filtered quantization of the Poisson algebra $\C[X]$, while $\C[\tldX]^{\bbW}$ is
the universal family of Poisson deformation of $\C[X]$. 
Let $S = \Abar(\sfD^{ch}(\tldX)) \cap \C[\frg^*]$ be a Poisson-commutative subalgebra of 
$\Abar(\sfD^{ch}(\tldX))$. Then, by 
the universality of $\C[\tldX]^{\bbW}$ (\cite[Proposition 2.12]{Losev-2016}), we have 
a unique homomorphism $\C[\frg^*]^{\bbW} \longrightarrow S$ and a unique 
isomorphism of Poisson algebras 
$\C[\tldX]^{\bbW} \otimes_{\C[\frg^*]^{\bbW}} S \isoto \Abar(\sfD^{ch}(\tldX))$. 
By \refprop{prop:Zhu-quantization}, the Zhu algebra $A(\sfD^{ch}(\tldX))$ is a 
filtered quantization of $\Abar(\sfD^{ch}(\tldX))$ over $S$. Thus, by 
\cite[Proposition 3.5]{Losev-2016}, we have a unique isomorphism 
$\sfD(\tldX)^{\bbW} \otimes_{\C[\frg^*]^{\bbW}} S \isoto A(\sfD^{ch}(\tldX))$.
Since we have the inclusions 
$\C[\tldX]^{\bbW} \subset \Abar(\sfD^{ch}(\tldX)) \subset \C[\tldX]$, the above 
homomorphisms are compatible with $\C[\frg^*]^{\bbW} \hookrightarrow S \hookrightarrow \C[\frg^*]$.
Thus we have the following proposition.

\begin{proposition}
 \label{thm:Zhu-alg}
 The Zhu algebra $A(\sfD^{ch}(\tldX))$ of the hypertoric vertex algebra 
 is a subalgebra of the quantized hypertoric algebra $\sfD(\tldX)$ which 
 contains its $\bbW$-invariant subalgebra $\sfD(\tldX)^{\bbW}$.
\end{proposition}

\footnotesize{
T.K.: Division of Mathematics, Faculty of Pure and Applied Sciences, University of Tsukuba,
Tsukuba, Ibaraki 305-8571, JAPAN. \\
{\em E-mail address}: \texttt{kuwabara@math.tsukuba.ac.jp}
}

\begin{thebibliography}{AKM}
\bibitem[A1]{Arakawa1}
T.~Arakawa,
{\it Vanishing of cohomology associated to quantized Drinfeld-Sokolov reduction},
Int. Math. Res. Not. {\bf 15} (2004), 729-767.

\bibitem[A2]{Arakawa2}
T.~Arakawa,
{\it Representation theory of $\scW$-algebras},
Invent. Math. {\bf 169} (2007), 219-320.

 \bibitem[AKM]{AKM}
T.~Arakawa, T.~Kuwabara and F.~Malikov,
{\it Localization of Affine W-algebras}, Comm.~Math.~Phys. {\bf 335} (2015), no.~1,
143-182.

\bibitem[BD]{Bielawski-Dancer}
R.~Bielawski and A.~S.~Dancer,
{\it The geometry and topology of toric hyperk\"{a}hler manifolds}, Comm. Anal. Geom.
{\bf 8}:4 (2000), 727-760.

\bibitem[BeKa]{Bezrukavnikov-Kaledin}
R.~Bezrukavnikov and D.~Kaledin,
{\it Fedosov quantization in the algebraic context}, Mosc. Math. J. {\bf 4} (2004),
559-592.

\bibitem[BeKu]{Bellamy-K}
G.~Bellamy and T.~Kuwabara,
{\it On deformation quantizations of hypertoric varieties},
Pacific J. Math. {\bf 260} (2012), no.~1, 89-127.

\bibitem[BLPW]{BLPW}
T.~Branden, A.~Licata, N.~Proudfoot, B.~Webster,
{\it Hypertoric category $\calO$}, Adv.~Math. {\bf 231} (2012), no. 3-4, 1487-1545.

\bibitem[FF1]{FF-Wakimoto}
B.~Feigin and E.~Frenkel,
{\it Affine Kac-Moody Algebras and Semi-infinite Flag Manifolds},
Comm. Math. Phys. {\bf 128} (1990), 161-189.

\bibitem[FF2]{Feigin-Frenkel}
B.~Feigin and E.~Frenkel,
{\it Affine Kac-Moody algebras at the critical level and Gel'fand-Diki\v{i} algebras},
Infinite Analysis, Part A, B (Kyoto, 1991), Adv. Ser. Math. Phys., vol. 16 (1992),
197-215.

\bibitem[FKW]{Frenkel-Kac-Wakimoto}
E.~Frenkel, V.~Kac and M.~Wakimoto,
{\it Characters and fusion rules for $\calW$-algebras via quantized Drinfel'd-Sokolov reduction},
Comm.~Math.~Phys. {\bf 147} (1992), no.~2, 295-328.

\bibitem[HS]{Hasel-Sturmfels}
T.~Hausel and B.~Sturmfels.
{\it Toric hyper{K}\"ahler varieties},
Doc. Math. {\bf 7} (2002), electronic, 495-534.

\bibitem[K]{BRST}
T.~Kuwabara, 
{\it BRST cohomologies for symplectic reflection algebras and quantizations of hypertoric 
varieties},
Transformation Groups {\bf 20} (2015), no.~2, 437-461.

\bibitem[KV]{Kaledin-Verbitsky}
D.~Kaledin and M. Verbitsky,
{\it Period map for non-compact holomorphically symplectic manifolds}, Geom. Funct. Anal.
{\bf 12} (2002), 1265-1295.

\bibitem[L1]{Losev}
I.~Losev,
{\it Isomorphisms of quantizations via quantization of resolutions},
Adv. Math. {\bf 231} (2012), 1216-1270.

\bibitem[L2]{Losev-2016}
I.~Losev,
{\it Deformations of symplectic singularities and orbit method for semisimple Lie algebras},
arXiv preprint (2016), \texttt{arXiv:1605.00592v1}.

\bibitem[MV]{Musson-Van-der-Burgh}
I.~M.~Musson and M.~Van~den Bergh,
{\it Invariants under tori of rings of differential operators and related
  topics},
Mem. Amer. Math. Soc. {\bfseries 136} (1998), no.~650, viii+85.

\bibitem[Wa]{Wakimoto}
M.~Wakimoto,
{\it Fock representation of affine Lie algebra $A_1^{(1)}$}, 
Comm. Math. Phys. {\bf 104} (1986), 605-609.

\bibitem[We]{Weibel}
C.~Weibel, {\it An introduction to homological algebra},
Cambridge studies in adv. math. {\bf 38}, Cambridge University Press, 1994.

\bibitem[Z]{Zhu}
Yongchang Zhu, 
{\it Modular invariance of characters of vertex operator algebras}, 
J.~AMS {\bf 9} (1996), num.~1, 237-302.

\end{thebibliography}
\end{document}